\documentstyle[twoside,12pt]{article}

\def\C          {{\cal C}}

\def\dim        {{\rm dim}}

\def\Char       {{\rm Char}}
\def\Spec       {{\rm Spec}}
\def\Supp       {{\rm Supp}} 
\def\Coker      {{\rm Coker}}
\def\Hom        {{\rm Hom}}
\def\Ker        {{\rm Ker}}

\newtheorem{prop}{Proposition}[subsection]
\newtheorem{dfn}[prop]{Definition}   
\newtheorem{theo}[prop]{Theorem}

\newtheorem{coro}[prop]{Corollary} 
\newtheorem{sublem}{Lemma}[subsection]

\title{\bf  Characteristic varieties of algebraic curves} 
\author{
{\sc {A.Libgober}}\\
{\small \it Department of Mathematics,
University of Illinois at Chicago, 
Chicago, Illinois, USA}}
\begin{document}

\maketitle 
\section*{0.Introduction.}

A procedure for calculation  of fundamental groups for the complements to 
algebraic curves in complex projective plane 
was found by Zariski and van Kampen ([Z],[vK]). 
Their methods yielded several important calculations and 
results on the fundamental groups of the complements 
(cf. for example [L2] for references). However,  
only limited information was obtained about their algebraic structure 
or what actually affects the complexity of these fundamental groups.
This paper is a result of attempts to find alternative 
ways for calculating the 
fundamental groups of the complements or at least some
invariants of these groups. The invariants of the fundamental groups,  
which we consider here, are certain subvarieties of complex tori 
${{\bf C}^*}^r$. They were called  characteristic varieties in [L3].
These subvarieties are unions of translated subtori, 
as follows from recent work [Ar].
We calculate these subtori in terms of local 
type of singularities and dimensions of linear systems which we 
attach to the configuration of singularities of the curve. 
\par These characteristic varieties  can be defined
as follows. Let ${\cal C}=\cup_{1 \le i \le r} C_i$ be an algebraic
curve in ${\bf C}^2$ and $\pi_1=\pi_1 ({\bf C}^2-{\cal C})$ be the
fundamental group of its complement. Then $\pi_1 /\pi_1'$ is isomorphic
to ${\bf Z}^r$ and acts on  $\pi_1'/\pi_1''$ by conjugation. This makes
$\pi_1'/\pi_1''$ into a module over the group ring of $\pi_1/\pi_1'$.
The latter is just the ring of Laurent polynomials 
${\bf Z}[t_1,t_1^{-1},...,
t_r,t_r^{-1}]$. After tensoring with ${\bf C}$, we obtain a 
${\bf C}[\pi_1/\pi_1']$-module
$\pi_1'/\pi_1'' \otimes {\bf C}$. The support of its 
$i-th$ exterior power is a subvariety of 
the torus ${\rm Spec} {\bf Z}[t_1,t_1^{-1},...,t_r,t_r^{-1}]$ 
called the i-th characteristic variety of $\cal C$: 
${\rm Char}_i({\cal C})$ (cf. [L3]). 
\par This invariant of the fundamental group can be used to calculate 
the homology of abelian covers of ${\bf C}^2-{\cal C}$ (cf. [L3]) 
and covers of ${\bf P}^2$ branched over the projective 
closure of ${\cal C}$
(cf. [Sa]). The above construction of characteristic varieties can be, 
of course,  
carried out for any topological space with $H_1={\bf Z}^r$
and several known results can be recasted using them. 
For example, the modules
$\pi_1'/\pi_1''$ were   
widely studied in the context of the complements to links in spheres
(cf. e.g. [Hil] and references there; for the case of algebraic
links cf. [Sab]). In this case,  
$\Char_1$ is the set of zeros of the multi-variable Alexander
polynomial. We shall see, however, that the characteristic varieties of 
algebraic curves rarely have codimension equal to one
and hence cannot be described using single polynomial. 
The varieties $\Char_i({\cal C})$
coincide with the cohomology support loci for local systems of rank 1 
considered in [Ar]. The homology of the Milnor fiber of the 
function obtained by homogenizing a defining equation of $\cal C$ 
(i.e. the Milnor fiber of the cone over the projective closure of $\cal
C$) can be found from the characteristic varieties of the latter. 
These Milnor fibers earlier were considered in the case when $\cal C$ is an  
 arrangement  i.e. when all components of ${\cal C}$ are lines (cf. [CS]).
\par For an irreducible $\cal C$  
the characteristic varieties  are 
subsets of  ${\bf C}^*$ i.e. collections of complex numbers. Those 
are the roots of the Alexander polynomial of $\cal C$
(cf. [L1]) and the results of this paper are equivalent to 
the results of [L2].
\par Our calculation of $\Char_i$ based on the following observations.
Firstly, by  Arapura's  theorem the characteristic varieties are translated 
tori and hence can be described by simple discrete data.
Secondly,  Sakuma's formula (cf. (1.3.2.2)), relating the homology of 
abelian branched
over $\cal C$ cover to the characteristic varieties of $\cal C$,
can be used to calculate such data for essential components (cf. 1.4.3)) 
of characteristic varieties completely
from the information about the homology of {\it all} abelian covers
with the branching locus $\cal C$. Thirdly, these abelian
covers can be realized as complete intersections and one can use 
the theory of adjoints (cf. (1.5)) to calculate the homology of these
covers (generalizing the calculations in the case of hypersurfaces cf. 
[Z], [L2], [L5];
complete intersections were used by Ishida (cf. [I], also [Zu]) in similar
context for calculations in the case of abelian covers of ${\bf P}^2$
branched along some arrangements).
We associate with each singular point a collection of polytopes in 
the unit cube $\cal U$ in ${\bf R}^r$ union of which is $\cal U$ and 
call them the local polytopes of quasiadjunction (cf. sect. 2.4.1).
Moreover, every such polytope defines the 
ideal in  the local ring of the singular point. The collection
of local polytopes defines new partition of the unit cube which is a 
refinement of partitions corresponding to all singular points of $\cal
C$ and reflecting the global information about singularities of $\cal
C$ (cf. sect. 2.6). 
We call the polytopes of this partition the global polytopes of
quasiadjunction. In the set of faces of global polytopes of 
quasiadjnction we single out a subset of contributing ones.
To each contributing face $\delta$ corresponds the linear system 
$H^0 ({\bf P}^2,{\cal J}_{\delta}(deg \ {\cal C}-3-l(\delta))$ 
where the ideal sheaf ${\cal J}_{\delta} \subset {{\cal O}}_{{\bf P}^2}$ 
and the integer $l(\delta)$ are  
determined by the face $\delta$. The components of $\Char_i{\cal C}$
correspond to contributing faces of global polytopes of quasiadjunction 
for which 
$dim  H^1 ({\bf P}^2,{\cal J}_{\delta}(deg \ {\cal C}-3-l(\delta)))=i$. 
The main result of the paper is the theorem 3.1 where the equations 
for translated tori are given explicitly in terms of corresponding faces of 
polytopes of quasiadjunction (cf. sect. 3).
\par The procedure for calculating the characteristic varieties,
though involving  possibly large calculations, is entirely algorithmic. 
For example, 
suppose that the curve $\cal C$ is an arrangement of $r$ lines.
(cf. sect. 3.3 for several
examples of explicite calculations for 
such arrangements including Ceva's and Hesse's arrangements).
Then any component of the characteristic variety $\Char_i({\cal C})$,
having positive dimension either belongs to
a component of characteristic variety of a sub-arrangement 
(i.e. is inessential and can be found by applying 
this algorithm to a sub-arrangement) or
is a connected component of a subgroup  corresponding to a collection 
${\cal S}$ of vertices in arrangement having  
multiplicity greater than $2$ (more detailed calculation allows 
to pick the component as well, cf. th.3.1). 
A collection ${\cal S}$ yields a component of characteristic 
variety if it satisfies the following conditions.
\par a) Certain system of linear homogeneous equations attached to 
$\cal S$ has a non zero solution. This system constructed as follows:
The unknowns $x_i$   
are in one to one correspondence with $r({\cal S})$ lines
of arrangement containing points from $\cal S$.   The equations are in 
one to one correspondence with the elements of ${\cal S}$. 
Left hand side of each equation is the sum of variables corresponding to
the lines through a point of ${\cal S}$ and the right hand side is 
a positive integer (not exceeding the number of lines in the
arramgements). 
\par b) The set  of solutions of the system from a) belongs to 
a hyperplane $\sum_{i=1}^{r({\cal S})} x_i=l({\cal S})$.
\par c) Let ${\cal I}_{\cal S}$ be the ideal sheaf with 
 $\Supp {\cal O}_{{\bf P}^2}/{\cal I}_{\cal S}=\cal S$ 
which stalk at $P \in {\cal S}$ 
is ${\cal M}_P^{m-2-\rho(P)}$ where ${\cal M}_P$ is 
the maximal ideal of $P$,  
$m$ is the number of lines in ${\cal S}$ containing $P$
and $\rho(P)$ is the right hand side of the equaition in the
system from a) corresponding to $P$. 
Then 
$dim H^1({\bf P}^2, {{\cal I}}_{\cal S} (r-3-l({\cal S}))) =i \ne 0$. 
\par Moreover if a)-c) are satisfied then 
the corresponding  to ${\cal S}$ subgroup is 
the set of solutions of 
$exp(L_P)=1$ where $L_P=\rho(P)$ are the equations 
of the system mentioned in a). Selection of particular connected 
component follows from more technical description in theorem 3.1
\par This algorithm yields complete description of characteristic tori
of positive dimension (and essential torsion points). 
An interesting problem which left
unanswered here is the realization problem: which tori (or collections
of tori) can appear as characteristic tori of an algebraic curve 
 with fixed degrees of components and given local type of
singularities. Some information in this direction however is provided in
section 4. Finally in the case of line arrangements (or
equivalently the case of fundamental groups of the complements
to arbitrary arrangements) the characteristic
varieties give new  sufficient conditions
(resonance conditions) for Aomoto complex of an arrangement to be 
quasi-isomorphic to the corresponding twisted DeRham complex  
(in many situations less restrictive than previously used (cf. [ESV]). 
We describe in a new way the space of ``resonant'' 
Aomoto complexes on given arrangement 
i.e. those with the cohomology different from the cohomology of 
generic Aomoto complexes 
(th. 5.4.1; this space was considered in [F]). Vice versa, this relation
between the space of resonant Aomoto complexes and characteristic
varieties shows that components of 
characteristic varieties which are subgroups of the group of 
characters are {\it combinatorial}
invariants of arrangements. Moreover, Aomoto complexes  provide 
another algorithm
for calculating these components of characteristic varieties
 of the fundamental groups of the complements to arrangements.
\par The announcement of these results is presented in [L6]. 
This work was supported by NSF grants DMS-9803623, DMS-9872025 
and Mittag Leffler Institute. I am also grateful to S.Yuzvinsky for 
very interesting correspondence regarding the material in section 5 and 
to J.Cogolludo for useful discussions of the manuscript.
\section{Preliminaries}
\subsection{Setting.} 
\smallskip Let $\bar {\cal C}$ $=\cup \bar C_i (i=1,...r)$ be a 
reduced algebraic
curve in ${\bf P}^2$ where $\bar C_i (i=1,...r)$ are the irreducible 
components of $\bar {\cal C}$.
We shall denote by $d_i$ the degree of the component $\bar C_i$. 
Let $L_{\infty}$ be a line in ${\bf P}^2$ which we shall view as the 
line at infinity. 
We shall be concerned with the fundamental groups of the 
complements to $\bar {\cal C}$ in ${\bf P}^2$ and in 
${\bf C}^2={\bf P}^2-L_{\infty}$. Let ${\cal C}=\cup C_i$ 
be the affine portion 
of $\bar {\cal C}$.
The homology of these complements are the following (cf. [L1]):
$$
H_1({\bf C}^2-{\cal C},{\bf Z})={\bf Z}^r, \ \ 
 H_1({\bf P}^2-\bar {\cal C},{\bf Z})={\bf Z}^r/(d_1,....,d_r) \eqno (1.1.1)
$$
Generators of these homology groups are represented by the
 classes of the 
loops $\gamma_i$, 
each of which is the boundary of a small 2-disk intersecting $C_i$ (resp. 
$\bar C_i$)
transversally at a non singular point.
\par For the fundamental groups we have the exact sequence:
 $$\pi_1({\bf C}^2-{\cal C}) \rightarrow \pi_1 ({\bf P}^2-\bar {\cal C}) 
\rightarrow 1 \eqno (1.1.2)$$
If the line $L_{\infty}$ is transversal to $\bar {\cal C}$, 
then the kernel of the 
surjection (1.1.2) is isomorphic to $\bf Z$ and belongs to the center of 
$\pi_1({\bf C}^2-{\cal C})$ (cf. [L4]). In general, the fundamental group
of the affine portion of the complement to $\bar \C$ in ${\bf P}^2$
depends on position of $L_{\infty}$ relative to $\bar \C$. Throughout the 
paper we assume that $L_{\infty}$ is transversal to $\bar \C$.
\medskip 
\subsection{Characteristic varieties of algebraic curves.}
\smallskip
\subsubsection{} 
Let $R$ be a commutative 
Noetherian ring and $M$ be a finitely generated 
$R$-module. Let  $\Phi: R^m \rightarrow R^n$ be such that
$M=\Coker \Phi$.  Recall that the  $k$-th Fitting ideal of $M$ is 
the ideal generated by $(n-k+1) \times (n-k+1)$ minors of the matrix 
of $\Phi$ (clearly depending only on $M$ rather than on $\Phi$).
The $k$-th characteristic variety $M$ is the reduced sub-scheme of 
$\Spec R$ defined by $F_k(M)$. 
\par If $R={\bf C}[H]$ where $H$ is a free abelian group then $R$ can be 
identified with the ring of Laurent polynomials and $\Spec R$ is a complex
torus. In particular each $k$-th characteristic variety of an $R$-module 
is a subvariety $V_k(M)$ of $({{\bf C}^*})^{{\rm {rk}} H} $.
\par If $Ann \wedge^k M \subset R$ is the annihilator of the 
$k$-th exterior power of $M$ then (cf. [BE], Cor.1.3):
 $(Ann \wedge^k M)^t \subseteq F_k(M) \subseteq Ann \wedge^k M$ for some
integer $t$. In particular, if $\Supp(M) \subset \Spec(R)$ is the set of prime
ideals in $R$ containing $Ann(M)$ (alternatively $\{\wp \in \Spec R 
\vert M \otimes R/\wp R \ne 0\}$, cf. [Se] p.3), then  
$\Supp (\wedge^k M)=\Supp(R/F_k(M))$ is the $k$-th characteristic
 variety of $M$.  
\par Note the following:    
\smallskip
\begin{sublem}  Let $0 \rightarrow M' \rightarrow M
\rightarrow M'' \rightarrow 0$. 
Then  $V_1(M)=V_1(M') \cup V_1(M'')$ and for $k \ge 2$: 
$V_k(M'') \subset V_k(M) \subset V_k(M'') \cup V_{k-1}(M'') \cap
V_1(M')$. 
\end{sublem}
\smallskip 
\par The first equality is Prop. 4(a) in [Se]. The second 
follows from the first and the exact sequence: 
$\Lambda^{k-1}M'' \otimes  M' \rightarrow \Lambda^k(M) \rightarrow 
\Lambda^k(M'') \rightarrow 0$, since $\Supp(A \otimes B)=
\Supp(A) \cap \Supp (B)$ for any $R$-modules of finite type
(cf. [Se], Prop. 4(c)).
\smallskip \par \noindent 
\subsubsection{} 
Let $G$ be a finitely generated, finitely 
presented group such that $H_1(G,{\bf Z})=G/G'={\bf Z}^r$
(for example $G=\pi_1({\bf C}^2-{\cal C})$ where
${\cal C}=\cup C_i$ is a plane curve as in 1.1; another class of 
examples which was studied in detail is given by link groups, cf. [Hil]).
Then  $G'/G'' \otimes {\bf C}$ can be viewed as 
$H_1(\tilde X,{\bf C})$ where $X$ is a topological space 
with $\pi_1(X)=G$ and $\tilde X$ is the universal abelian cover of $X$. 
The group $G/G'=H_1(X,{\bf Z})$ acts as the group of deck
transformations on $\tilde X$ and hence $G'/G'' \otimes {\bf C}$
has a structure of a ${\bf C}[G/G']$-module.  
We shall denote the i-th characteristic variety  of this module 
as $V_i(G)$ (or $V_i({\cal C})$ if $G=\pi_1({\bf C}^2-{\cal C})$)
and call it the i-th characteristic variety of $G$ (resp. $\cal C$).
The {\it depth}
of a component $V$ is the integer 
$i={\rm max} \{j \vert V \subset V_j(G)\}$.
We shall see below that if a component has depth $i$ and 
dimension $\varrho>0$ and contains identity, then 
$i=\varrho-1$, cf. footnote in 1.4.2.
\smallskip \par \noindent 
{\rm 1.2.2.1.} $\  $
\smallskip \par \noindent
If $G=F_r$ is a free group on $r$-generators then $G'/G''
=H_1(\widetilde {\bigvee_r S^1},{\bf Z})$, where $\widetilde {\bigvee_r S^1}$
is the universal abelian cover of the wedge of $r$ circles.  
It fits into the exact sequence: 
\begin{displaymath}
0 \rightarrow 
H_1(\widetilde {\bigvee_r S^1},
{\bf C}) \rightarrow {\bf C}[{\bf Z}^r]^r \rightarrow I \rightarrow 0
\end{displaymath}  
with $I$ denoting the augmentation ideal of the group ring of ${\bf Z}^r$.
(As an universal abelian cover of $\widetilde {\bigvee_r S^1}$ one can take
the subset of ${\bf R}^r$ of points having at least $r-1$ integer
coordinates with the action of ${\bf Z}^r$ given by translations; unit 
vectors of the standard basis provide identification of 1-chains 
with ${\bf C}[{\bf Z}^r]^r$ while the module of 0-chains is identified 
with ${\bf C}[{\bf Z}^r]$). This shows that 
$H_1(\widetilde{\bigvee_r S^1},{\bf C})$ is cokernel of the 
map $\Lambda^{r \choose 3}{\bf C}[{\bf Z}^r]^r \rightarrow 
\Lambda^{r \choose 2} {\bf C}[{\bf Z}^r]^r$ in the Koszul resolution
corresponding to the $(x_1-1),...,(x_r-1)$.  
This implies that 
$V_i(F_r)={{\bf C}^*}^r$ for $0 < i \le r-1$ and $V_i(F_r)=(1,...,1)$
for $r \le i \le {r \choose 2}$ (cf. also (1.4.1) below). 
\par If $G$ is the fundamental group of a link in a 3-sphere $S^3$ 
with $r$ components,  
then the first determinantal ideal is generated by 
$(t_1-1),...,(t_r-1)$ and certain principal ideal. A generator
$\Delta(t_1,...,t_r)$ of the latter is called the Alexander polynomial.
Alexander polynomial satisfies $\Delta(1,...,1)=0$ 
and hence $V_1(G)$ is the hypersurface $\Delta(t_1,...,t_r)=0$. 
Extensive calculations of the Alexander polynomials of links can be 
found in [SW].
In particular, if $G$ is the fundamental group of the complement to
the Hopf link 
in $S^3$ with $r$ components then $V_1(G)$ is the set of zeros of 
$t_1\cdot t_2 \cdot \cdot \cdot t_r-1$.
Moreover $V_1(G)=...=V_{r-1}(G)$ (cf. [L4],p.165). 
From a presentation of $G$ using Fox calculus one can calculate a presentation
of $\pi_1'/\pi_1'' \otimes {\bf C}$ as a ${\bf C}([H])$-module and hence
the characteristic varieties of $\pi_1$ (cf. [Hi],[CS] for examples of such 
calculations).
\smallskip \par \noindent 
\subsubsection{}
Let  $T(L_{\infty})$
be a small tubular neighborhood of $L_{\infty}$  in ${\bf P}^2$. If  
$\partial T(L_{\infty})$ is its boundary then 
$\pi_1 (\partial T(L_{\infty})-\partial T(L_{\infty}) \cap {\cal C})
\rightarrow \pi_1({\bf P}^2-L_{\infty} \cup {\bar \C})$ is 
a surjection. Lemma 1.2.1 implies that the characteristic 
variety of $\cal C$ is a subset of the torus $t_1^{d_1} \cdot \cdot \cdot
t_r^{d_r}=1$ since ${\cal C} \cap \partial T(L_{\infty}) \subset
 \partial T(L_{\infty})=S^3$ is the Hopf link with $d_1+...+d_r$
components with $d_i$ components of the link belonging to  
$C_i$ and hence corresponding to $t_i$ for each $1 \le i \le r$ (cf. [L3]). 
In fact  
the characteristic varieties of an affine 
curve can be determined from the projectivization
 as follows. 
\par It is a corollary of (1.1.1) that ${\bf T}_a=
\Spec{\bf C}[H_1({\bf C}^2-{\cal C})]$ is the torus of dimension $r$ and
that ${\bf T}_p=\Spec {\bf C}[H_1({\bf P}^2-{\bar \C})]$ is the sub-scheme
of zeros of $t_1^{d_1} \cdot \cdot \cdot t_r^{d_r}-1$ in ${\bf T}_a$.
We denote by $E: {\bf T}_p \rightarrow {\bf T}_a$ the corresponding embedding.
On the other hand, the construction of (1.2.1) and (1.2.2) yields 
subvarieties  
$V_i ({\bar \C})_{p}$ in ${\bf T}_p$.
\smallskip 
\par \noindent
\addtocounter{prop}{2} 
\begin{prop} The characteristic variety 
of projective and affine curves satisfy: 
$$V_i({\cal C})=E(V_i({\bar \C})_p) \eqno (1.2.3.1)$$
\end{prop} 
\par \noindent
{\bf Proof}. It follows from the isomorphism: 
     $$\pi_1'({\bf P}^2-{\bar \C})/\pi_1''({\bf P}^2-{\bar \C})=
         \pi_1'({\bf C}^2-{\cal C})/\pi_1''({\bf C}^2-{\cal C})
 \eqno   (1.2.3.2)$$
equivariant with respect to the action of $H_1({\bf C}^2-{\cal C})$.
This isomorphism is a consequence of (1.1.2) because in the latter
the left map 
induces isomorphism on commutators. Indeed, the kernel of 
surjection (1.1.2) is isomorphic to ${\bf Z}$ (cf. (1.1)) and 
does not intersect $\pi_1'({\bf C}^2-{\cal C})$ 
because it injects into $H_1({\bf C}^2-{\cal C})$ 
(cf. also [L4]).

\smallskip \noindent 
\subsection{Abelian covers}
\subsubsection{}
Let $m_1,...,m_r$ be positive integers and 
$h_{m_1,..,m_r}: H_1({\bf C}^2-{\cal C},{\bf Z}) \rightarrow 
{\bf Z}/m_1 {\bf Z} \oplus ...\oplus {\bf Z}/m_r {\bf Z}$ be 
the surjection $\gamma_i \rightarrow \gamma_i \ mod \ m_i$.
The  kernel of the homomorphism $\pi_1 ({\bf C}^2-{\cal C}) 
\rightarrow {\bf Z}/m_1 {\bf Z} \oplus ...\oplus {\bf Z}/m_r {\bf Z}$,
which is the composition of the abelianization $ab: \pi_1({\bf C}^2-{\cal C}) 
\rightarrow H_1({\bf C}^2-{\cal C})$ and $h_{m_1,...,m_r}$, 
defines an unbranched cover of ${\bf C}^2-{\cal C}$. We shall denote
it as $\widetilde {({\bf C}^2-{\cal C})_{m_1,...,m_r}}$.
This is a quasi-projective algebraic 
variety defining a birational class of projective surfaces 
$\overline {\widetilde {({\bf C}^2-{\cal C})_{m_1,...,m_r}}} $. Birational 
invariants of surfaces in this class (in particular the first 
Betti number of a non singular model) depend only on ${\cal C}$
and the homomorphism $h_{m_1,...,m_r}$.
\par If $h_{m_1,...,m_r}(d_1\gamma_1+...d_r\gamma_r)=0$,
then the corresponding branched covering of ${\bf C}^2$
 is a restriction of the covering of ${\bf P}^2$ unbranched over the 
line at infinity. It can be easily checked that the first Betti numbers 
of those two branched coverings are the same,  since we 
assume (cf. (1.1))
 that the line at infinity is transversal to $\cal C$ 
(cf. [L1]).
\par A model (singular, in general) for a surface birational to   
$\overline  {\widetilde {({\bf C}^2-{\cal C})_{m_1,...,m_r}}}$
can be constructed as follows.
Let $f_i(u,x,y)=0$ be an equation of the component $C_i$ ($i=1,...,r$). 
Let $V_{m_1,...,m_r}$ be a complete intersection on ${\bf P}^{r+2}$ 
(coordinates of which we shall denote $z_1,..,z_r,u,x,y$)
given by the equations
 $$z_1^{m_1}=u^{m_1-d_1}f_1(u,x,y),..., 
z_r^{m_r}=u^{m_r-d_r}f_r(u,x,y) \eqno (1.3.1.1)$$
Projection from the subspace given by $u=x=y=0$ onto the plane
$z_1=...=z_r=0$ (i.e. $(z_1,..,z_r,u,x,y) \rightarrow (u,x,y)$), when 
restricted on the preimage in $V_{m_1,...,m_r}$ of  ${\bf C}^2-{\cal C}$,
is unbranched cover of ${\bf C}^2-{\cal C}$ 
corresponding to $Ker (h_{m_1,...,m_r} \circ ab)$. 
\smallskip \par \noindent 
\subsubsection{}
The first Betti number of {\it unbranched} cover 
$ {\widetilde {({\bf C}^2-{\cal C})_{m_1,...,m_r}}}$ can be  
found in terms of the characteristic varieties of $\cal C$ as follows
(cf. [L3]). For $P \in {{\bf C}^*}^r$ let 
$f(P,{\cal C})=\{max \quad i \vert P \in V_i
({\cal C}) \}$. Then 
$$b_1({\widetilde {({\bf C}^2-{\cal C})_{m_1,...,m_r}}})=
r+\Sigma_{\omega_i^{m_i}=1, (\omega_{m_1},...,\omega_{m_r}) \ne
(1,...,1)} f((\omega_{m_1},..,\omega_{m_r}),{\cal C}) \eqno (1.3.2.1)$$ 
\par The first Betti number of a resolution of branched cover of ${\bf P}^2$
(i.e. of $V_{m_1,...,m_r}$) can be calculated using the characteristic
varieties of curves formed by components of $\cal C$ (cf. [Sa]). 
Let $\widetilde V_{m_1,..,m_r}$ be such a resolution. 
For a torsion point
of $\omega=(\omega_1,...,\omega_r), \omega_{i}^{m_i}=1$
in the torus  ${{\bf C}^*}^r$ let
${\cal C}_{\omega}=\cup_{i \vert \omega_i \ne 1} C_i$. 
Then the first Betti number of $\widetilde V_{m_1,...,m_r}$
equals:
$$\Sigma_{\omega} max \{ i \vert \omega \in Char_i (\C_{\omega}) \}
 \eqno (1.3.2.2)$$
More precisely, if $\chi_{\omega}$ is the character of 
$\pi_1({\bf C}^2-{\cal C})$ such that $\chi_{\omega}(\gamma_i)=\omega_i$ 
and for a character $\chi$ of the Galois group  
$Gal(\widetilde V_{m_1,...,m_r}/{\bf P}^2)$ we put: 
$$H_{1,\chi}(\widetilde V_{m_1,...,m_r})=\{ x \in H_1(\widetilde
 V_{m_1,...,m_r}) \vert g(x)=\chi (g)\cdot x, \forall g \in 
 Gal(\widetilde V_{m_1,...,m_r}/{\bf P}^2)  \}
 \eqno (1.3.2.3)$$ then 
 $$dim H_{1,\chi_{\omega}}=max \{ i \vert \omega \in Char_i (C_{\omega})\}
 \eqno (1.3.2.4)$$ 
\smallskip \par \noindent
\subsubsection{A bound on the growth of Betti number.}
\smallskip \par \noindent
\addtocounter{prop}{2}
\begin{prop} Let $ b_1(\bar {\cal C},n)$ (resp.
 $ {b_1({\cal C},n)}$)  be the first Betti
number of the cover of ${\bf P}^2$  (resp. ${\bf C}^2-{\cal C}$) 
branched over $L_{\infty} \cup {\bar \C}$ (resp. unbranched) 
and corresponding to the surjection $h_{n,...,n}: \pi_1({\bf P}^2-
L_{\infty} \cup {\bar \C}) \rightarrow ({\bf Z}/n{\bf Z})^r$ (given by 
evaluation modulo $n$ of the linking numbers of loops with the 
components of $\cal C$ modulo $n$). Then ${b_1({\bar {\cal C},n)}} \le
\bar C_1 \cdot n^{r-1}$. (resp. $ {b_1({\cal C},n)} \le C_1
\cdot n^{r-1}$) for some constants $C_1,\bar C_1$
independent of $n$.
\end{prop}
\smallskip \par \noindent
{\bf Proof}.
This follows from the Sakuma's formula (1.3.2.2)
(resp. (1.3.2.1)) and the obvious remark that the number of 
$n$-torsion points on a torus of dimension $l$ grows as $n^l$ since  
$dim (Char_i(\pi_1({\bf C}^2-{\cal C})'/\pi_1({\bf C}^2-{\cal C})'') \le
r-1$ by 1.2.3.
\smallskip \par \noindent 
\subsubsection{Characteristic varieties and the homology of Milnor
fibers.} 
\smallskip \par The polynomial $f_1 (u,x,y)\cdot \cdot \cdot f_r(u,x,y)$ 
(which set of zeros in ${\bf P}^2$
is $\bar \C$) defines a cone in ${\bf C}^3$ having a non isolated 
singularity, provided $\cal C$ is singular. The Milnor fiber $M_c$ of this 
singularity (cf. [CS] in the case when $degf_i=1, \forall i$)
is diffeomorphic to an affine hypersurface given by the
equation: $f_1 \cdot \cdot \cdot f_r=c, c \ne 0$.
Quotient of the latter by the 
action of the cyclic  group ${\bf Z}/d{\bf Z}$ ($d=\Sigma_i d_i, d_i=deg f_i$)
acting via $(u,x,y) \rightarrow (\omega_d u, \omega_d x, \omega_d y),
\omega_d^d=1$ is ${\bf P}^2-{\bar \C}$. In other words, the Milnor fiber is the
cyclic cover $p: M_c \rightarrow {\bf P}^2-{\bar \C}$ 
corresponding to the homomorphism sending
$\gamma_i \rightarrow 1 \quad mod \quad d$. The exact sequence of 
the pair $(M_c, p^{-1}({\bf P}^2-{\bar \C} \cup L_{\infty}))$ shows that 
$rk H_1(M_c)=rk H_1(p^{-1}({\bf P }^2-{\bar \C} \cup L_{\infty})-1$,
since we assume that $\bar \C$ is transversal to $L_{\infty}$ 
(cf. [L1]). Hence it follows from
(1.3.2.1) that 
$$rkH_1(M_c)= r-1+\Sigma_{i=1}^{d-1} f((\omega_{d}^i,..,\omega_d^i),{\cal C})
 \eqno (1.3.4.1)$$  
\smallskip \par \noindent
\subsection{Characteristic varieties and support loci for rank one local
systems}
\subsubsection{}
Let again $G$ be a group such that $G/G'={\bf Z}^r$. If $X$ is
a topological space with $\pi_1(X)=G$ then the local systems of rank one
on $X$ correspond to the points $\Hom (G, {\bf C}^*)$ (cf. [St]). The latter 
has a natural identification with $H^1(X,{\bf C}^*)$.
Each $\gamma_i$, corresponding to a component $C_i$ of ${\cal C}$ 
(cf. 1.1), defines the homomorphism 
$t_i: \Hom(G,{\bf C}^*) \rightarrow {\bf C}^*$ given by $t_i(\chi)=
\chi(\gamma_i), \chi \in \Hom(G,{\bf C}^*)$. Therefore $t_i$'s  provide
 an identification of 
$\Hom(G,{\bf C}^*)$ with ${{\bf C}^*}^r$.
\par The homology groups $H_i(X,\rho)$ of $X$ with coefficients in a local
system  corresponding to a homomorphism
$\rho: \pi_1(X) \rightarrow H_1(X,{\bf Z}) \rightarrow {\bf C}^*$
 are the homology of the complex
$C_i(\tilde X) \otimes_{H_1(X,{\bf Z})} {\bf C}$ where $\bf C$ is
equipped with the structure of ${\bf Z}[H_1(X,{\bf Z})]$-module using
$\rho$. If $\rho \ne 1$ then 
 $$H_1(\tilde X,{\bf C}) \otimes_{{\bf C}[H_1(X,{\bf Z})]} {\bf C}= 
 H_1(X,\rho) \eqno (1.4.1.1)$$ This follows, for example, from the exact sequence of 
the low degree terms in the spectral sequence corresponding 
to the action of $H_1(X,{\bf Z})$ on the universal abelian cover 
$\tilde X$: $H_p(H_1(X,{\bf Z}),H_q(\tilde X)_{\rho}) \Rightarrow
H_{p+q}(X,\rho)$  (here $H_q(\tilde X)_{\rho}$ denotes 
the homology of the complex $C_i(\tilde X) \otimes_{\bf Z}{\bf C}$
with the action of $H_1(X,{\bf Z})$ given by $g(e \otimes \alpha)=
g \cdot e \otimes \rho(g^{-1})\alpha, g \in H_1(X,{\bf Z}), e \in C_i(\tilde
X),\alpha \in {\bf C}$ i.e. the  usual 
homology $H_q(\tilde X,{\bf C})$ with the action of $H_1(X,{\bf Z})$
changed by the character $\rho$ (cf. [CE],ch. XVI,th. 8.4).
This exact sequence is:
$$H_2(X,\rho) \rightarrow H_2(H_1(X,{\bf Z}),\rho) 
\rightarrow (H_1(\tilde X)_{\rho})_{H_1(X,{\bf Z})} \rightarrow
H_1(X,\rho) \rightarrow H_1(H_1(X,{\bf Z}),\rho) \rightarrow 0 $$
(cf. [CE], ch XVI, (4a)).
Since for $\rho \ne 1$ we have  $H_i(H_1(X,{\bf Z}),\rho)=0$, we obtain 
(1.4.1.1). For $\rho=1$, an argument similar to [L3], sect. 1 yields
that $dim H_1(\tilde X,{\bf C})\otimes_{{\bf C}[H_1(X,{\bf Z})]} 
{\bf C}$ is the dimension of the kernel of the map
$\cup_X: \Lambda^2 H^1(X,{\bf C}) \rightarrow H^2(X,{\bf C})$ 
given by the cup product. 
From the definition of Fitting ideals (cf. 1.2.1)
it follows that for $\rho \ne 1$ one has:
 $$V_i(X)=\{ \rho \in \Hom(G,{\bf C}^*) \vert H_1(X,\rho) \ge i \}
 \eqno (1.4.1.2)$$
and that $\rho=1$ belongs to $V_{dim \Ker \cup_X}$ (cf. [L3], Prop.1.1).
\par For example if $G=F_r$ then $dim H_0(F_r,\rho)$ is $0$, 
if $\rho$ is non trivial, and $1$ otherwise. Using $e(F_r,\rho)=r-1$
we obtain that $dim H^1(F_r,\rho)$ is $r-1$, if $\rho$ is non trivial, 
and otherwise is $r$. Since $dim Ker \cup_{F_r}={r \choose 2}$ 
we recover the description of the characteristic
varieties for $F_r$ mentioned in (1.2.2.1).  
\subsubsection{Structure of characteristic varieties.} 
\smallskip \par We will need
the following theorem of D.Arapura (cf. [Ar]) which generalizes the 
results of C.Simpson to quasi-projective case.
\par Let $\bar X$ be a K\"ahler manifold such that 
$H^1(\bar X,{\bf C})=0$,
$D$ a divisor with normal crossings and $X=\bar X-D$. Then, 
for each characteristic variety $V$, there exist a
finite number of torsion characters $\rho_i \in \Hom(G,{\bf C}^*)$, a
finite number of unitary characters $\rho_j'$ and surjective maps
onto (quasiprojective) curves $f_i: X \rightarrow C_i$ such that 
 $$V(X)=\bigcup_i \rho_if^*H^1(C_i,{\bf C}^*) \cup \bigcup \rho_j'
 \eqno (1.4.2.1)$$ 
\par A consequence of 1.4.2.1 for curves in ${\bf C}^2$ is that the
components of positive dimensions of their characteristic varieties are subtori
of ${{\bf C}^*}^r$ translated by points of finite order
\footnote{{though the paper [Ar] considers only the case of the first
characteristic variety (i.e. in terminology of [Ar] characters $\rho$
such that $dim H^1(\rho) \ge 1$), D.Arapura communicated to the author
that the statement is true for all $V_k$. Moreover it follows from his 
argument that the dimension of $V_k$, {\it containing the identity of 
the group of characters}, is $k+1$. Indeed by (1.4.2.1) for 
any local system in such irreducible component of positive dimension
of $V_k(X)$ there exist $L'$ on an appropriate curve $C$ and the map
 $f: X \rightarrow C$ such that $L=f^*L'$. Moreover it follows from 
Proposition 1.7 in [Ar] that for all but finitely many $L$ one has
 $H_1(X,L)=H_1(C,L')$. But $\pi_1(C)$ is free and if $k+1$ is the 
number of its generators and $L'$ is not trivial then 
$dim H_1(C,L')=k$ (cf. 1.2.2.1). }}.
\subsubsection{Essential for a given set of components tori.} 
\par By coordinate
torus (corresponding to components $C_{i_1},...
C_{i_s}$) we shall mean a subtorus in ${{\bf C}^*}^r$ given by 
$$t_{i_1}=...=t_{i_s}=1. \eqno (1.4.3.1)$$
The inclusion $I_{i_1,..,i_s}: 
{\bf C}^2-\cup_{i=1,..,r} C_i \rightarrow 
{\bf C}^2-\cup_{i \ne i_1,...,i_s}  C_i$ induces a 
surjective map  
$\tilde I_{i_1,..,i_s}: \pi_1 ({\bf C}^2-\cup_{i=1,...,r} C_i)
 \rightarrow \pi_1 ({\bf C}^2-\cup_{i \ne i_1,...,i_s} C_i$) with 
restriction $\tilde I'_{i_1,..,i_s}: \pi_1'({\bf C}^2-\cup_{i=1,...r}C_i) 
\rightarrow \pi_1'({\bf C}^2-\cup_{i \ne i_1,..i_s}C_i)$ which is 
also surjective. Indeed if $K=\Ker \pi_1({\bf C}^2-\cup_{i=1,,,r}C_i)
\rightarrow H_1({\bf C}^2-\cup_{i \ne i_1,...,i_s}C_i)$ then 
$K \rightarrow \pi_1'({\bf C}^2-\cup C_{i \ne i_1,..,i_s})$ 
is surjective. Since $K'$ is a normal closure of 
$\pi_1'({\bf C}^2-\cup_{i=1,...r} C_i)$ and loops trivial in 
$\pi_1({\bf C}^2-\cup_{i \ne i_1,...,i_s} C_i)$ (e.g. loops which 
consist of paths from the base point to a point in vicinity of 
$C_i, (i \ne i_1,..,i_s)$  and loops bounding small disk transversal to 
$C_i$) the surjectivity of $I_{i_1,..,i_s}'$ follows.
\par The latter gives rise to a surjective map of
 ${\bf C}[H_1( {\bf C}^2-\cup_{i=1,..,r} C_i)]$-modules:
\break ${\pi_1'/\pi_1''({\bf C}^2-\cup_{i=1,..,r} C_i) 
\rightarrow \pi_1'/\pi_1''({\bf C}^2-\cup_{i \ne i_1,...,i_s}  C_i)}$
which induces an injection of corresponding characteristic varieties:
 $$V_k({\bf C}^2-\cup_{i \ne i_1,...,i_s}  C_i)
 \rightarrow V_k ({\bf C}^2-\cup_{i=1,..,r} C_i) \eqno (1.4.3.2)$$
(cf. Lemma 1.2.1). A component of $V_k ({\bf C}^2-\cup_{i=1,..,r} C_i)$ 
which is an image of a component for 
some $i_1,...,i_s$ in (1.4.3.2) is called obtained via a pull back. 
A component of $V_k ({\bf C}^2-\cup_{i=1,..,r} C_i)$ is 
called essential if it isn't a pull back of component of 
a characteristic variety of a curve composed of irreducible components 
of $\cal C$.
\smallskip 
\par \noindent
\addtocounter{sublem}{2}
\begin{sublem} Let $V$ be a connected component of the 
characteristic variety $V_1$ of ${\cal C}$ having positive dimension 
and 
belonging to the coordinate torus $t_{i_1}=...=t_{i_s}=1$.
Then it is obtained via a pull back of a component of characteristic variety
for the union of components of ${\cup C_i}$
 ($i \ne i_1,...,i_s$).
\end{sublem} 
\smallskip \par \noindent
{\bf Proof}. According to Arapura's theorem (cf. [Ar],(1.4.2)) 
 component $V$
defines a map $f: {\bf C}^2- {\cal C} \rightarrow C$ for some 
quasiprojective curve $C$ such that for some local system 
$E \in \Char \pi_1({\bf C}^2 -\C)$ 
one has: $V=E\otimes 
f^* (\Char C)$ where $\Char C=
\Hom(\pi_1(C),{\bf C}^*)$.  
We claim that $f$ factors as follows: 
 $$\matrix {{\bf C}^2-{\cal C} & \buildrel I_{i_1,...,i_s} \over 
\longrightarrow & 
 {\bf C}^2-\cup_{i \ne i_1,...,i_s} C_i \cr
 & \searrow f & \downarrow  \tilde f  \cr
 & &  C} \eqno (1.4.3.3)$$
The lemma is a consequence of  existence of $\tilde f$. Indeed for almost 
all local systems $L$ on $C$ we have $H^1(E \otimes f^*(L))=
H^1(f_*(E) \otimes L)$ (cf. proof of Prop. 1.7 in [Ar]).
Moreover $H^1(f_*(E) \otimes L)=
H^1(\tilde f_* \circ (I_{i_1,..,i_s})_* (E) \otimes L)=
H^1((I_{i_1,..,i_s})_*E \otimes \tilde f^*L)$  and the latter has 
the same dimesnion for almost all $L$ again by the same argument from 
the proof of Prop. 1.7 in [Ar].
\par To show the existence of $\tilde f$, let $D=\bar C-C$ where $\bar C$ is a non singular compactification
of $C$. Since for $j=i_1,..,i_s$ we have $t_j=1 $ on a 
translate (i.e. a coset) of $f^*(\Hom (H_1(C),{\bf C}^*))$  and hence $t_j=1$
on the latter subgroup of $\Char(\pi_1({\bf C}^2-{\cal C}))$ we have  
for any $\chi \in \Char(H_1(C,{\bf Z}))$ and $j=i_1,..,i_s$ the following: 
$\chi(f_*(\gamma_j))=f^*(\chi)(\gamma_j)=t_j(f^*(\chi))=1$.
Equivalently $f_*(\gamma_j)=0$. Therefore $f_*(\gamma)=0$ 
in $H_1(C,{\bf Z})$ for 
 $\gamma \in H_1({\bf C}^2-{\cal C},{\bf Z})$ 
if and only if $\gamma$ belongs to the 
subgroup generated by $\gamma_{i_1},...,\gamma_{i_s}$.
Let us consider the pencil of curves on ${\bf P}^2$ formed by the 
fibres of $f$. $f$ extends to the map from the complement to the base locus 
of this pencil to $\bar C$ 
\footnote{{incidentally, since the resolution $\bar X$ 
of the base locus of this pencil is simply-connected one has 
$\bar C={\bf P}^1$}}. Preimage of $D$ 
in this extension is a union of components of $\C$ 
and we want to show that none of these components is $C_i$ with 
$i=i_1,...i_s$. 
But none of the components $C_i, i=i_1,...,i_s$ is taken by this 
extension into $D$ since otherwise $f_*(\gamma_i) \ne 0$ for the 
corresponding $\gamma_i$. Hence domain of this extension of $f$ 
contains all points of $C_{i_1} \cup ... \cup C_{i_s}$ 
not belonging to the remaining components of $\C$.
\par  Note that for isolated points of characteristic varieties
it can occur that 
$H^1({\bf C}^2-\cup_{i \ne i_1,..,i_s} C_i,L) \ne
H^1({\bf C}^2-C,I^*_{i_1,..,i_s}(L))$ as is shown by examples in [CS].
\subsection{Adjoints for complete intersections}
\subsubsection{} Let $F \subset {\bf P}^{n}$ be a surface which is a
complete intersection given by the equations:
 $$\bar F_1=...=\bar F_{n-2}=0 \eqno (1.5.1.1)$$ 
of degrees $d_1,...,d_{n-2}$ respectively.
Let (cf. [Ha] p.242)
$$\Omega_F={\cal E}xt^{n-2}({\cal O}_F,
{\Omega}^n_{{\bf P}^n})) $$ 
be the dualizing sheaf of $F$. From the latter and the Koszul resolution
$$0 \rightarrow {\cal O}_{{\bf P}^n}(-d_1-...-d_{n-2}) \rightarrow 
 .... \rightarrow  {\cal O}_{{\bf P}^n}(-d_1) \oplus ... \oplus 
{\cal O}_{{\bf P}^n}(-d_{n-2}) \rightarrow {\cal O}_{{\bf P}^n}
\rightarrow {\cal O}_F \rightarrow 0 \eqno (1.5.1.2)$$
it follows that one can identify $\Omega_F$ with 
${\cal O}_F(d_1+...+d_{n-2}-n-1)$.
\par Let $f: \tilde { F} \rightarrow  F$ be a resolution of 
singularities of $F$ and 
$\tau: f_*(\Omega_{\tilde { F}}) \rightarrow \Omega_{ F}$ be the
trace map (cf. [BL]).
 It identifies sections of $f_*(\Omega_{\tilde F})$ over an
open set with those meromorphic differentials on non singular part of this
open set in $F$  which when pulled back on a resolution $\tilde F$
admit a holomorphic extension over the exceptional set of $f$.
The adjoint ideal ${\cal A}'$ is the annihilator of the cokernel of 
$\tau$ i.e. 
$${\cal A}'={\cal H}om_{{\cal O}_F}(\Omega_F, f_*(\Omega_{\tilde F}))=
f_*(\Omega_{\tilde F})(-d_1-...-d_{n-2}+n+1) \eqno (1.5.1.3)$$
We define the sheaf of adjoint ideals on ${\bf P}^n$ as
 ${\cal A}=\pi^{-1}({\cal A}')$ (also denoted as $Adj F$) 
where $\pi$ is the most right map in (1.5.1.2).
The degeneration of Leray spectral sequence for $f$ (due to the 
Grauert-Riemenschneider vanishing theorem ([GR])) yields 
$$H^i(\tilde F,\Omega_{\tilde F})=H^i(F,f_*(\Omega_{\tilde F}))=$$ 
$$=H^i(F,{\cal A}'(d_1+...+d_{n-2}-n-1)=
H^i({\bf P}^n,{\cal A}(d_1+...+d_{n-2}-n-1))\eqno (1.5.1.4)$$
In particular the irregularity of $\tilde F$ i.e. 
$dim H^1(\tilde F,{\cal O})$ can be found as the difference between the
actual dimension $H^0({\bf P}^n,{\cal A}(d_1+...+d_{n-2}-n-1)$) 
and the ``expected'' dimension (i.e. $\chi({\cal A}(d_1+...d_{n-2}-n-1))$)
of the adjoints (since $H^i({\cal A}(d_1+...d_{n-2}-n-1))=0$ for 
$i \ge 2$).
\subsubsection{Local description of adjoint ideals.} Let 
$$F_1(w_1,....,w_n)=0,...,F_{n-2}(w_1,...,w_n)=0 \eqno (1.5.2.1)
$$ be a germ of a complete intersection of hypersurfaces
in ${\bf C}^n$ having an isolated singularity at the origin $O$. 
For any two pairs $1 \le i,j \le n, 
i \ne j$ and $1 \le k,l \le n, k \ne l$ we have up to sign:
$$ {{dw_i \wedge dw_j} \over {{\partial (F_1,...,F_{n-2})} \over 
{\partial (w_1,...,\hat w_i,..,\hat w_j,..,w_n)}}}=
 {{dw_k \wedge dw_l} \over {{\partial (F_1,...,F_{n-2})} \over
 {\partial (w_1,...,\hat w_k,...,\hat w_l,...,w_n)}}} \eqno (1.5.2.2)$$
Indeed the Cramer's rule for the solutions of the system of equations:
 $${{\partial F_k} \over {\partial w_1}} dw_1 \wedge dw_i +...+
 {{\partial F_{k}} \over {\partial w_n}} dw_n \wedge dw_i=0 \ \ 
(k=1,...,n-2)$$
when one views $dw_k \wedge dw_i (k=1,..,\hat i,.,n-1)$ as unknowns  
yields that up to sign: 
$$ dw_k \wedge dw_i={{{{\partial (F_1,...,F_{n-2}) \over
 {\partial (w_1,...,
\hat w_k,...,\hat w_i,...,w_n)}} dw_n \wedge dw_i} \over
 {{{\partial (F_1,...,F_{n-2})} \over
 {\partial (w_1,...,\hat w_i,...,w_{n-1})}}}}} \eqno (1.5.2.3)$$
(1.5.2.2) follows from this for any two pairs $(i,j),(k,l),i \ne j,k \ne l$.
\par Since $F_1=....=F_{n-2}=0$ is a complete intersection with isolated
singularity one of the Jacobians ${\partial (F_1,..,F_{n-2})}
 \over {\partial (w_1,..,\hat w_i,...,\hat w_j,...,w_n)}$ is non 
vanishing in a neighborhood of the singularity everywhere except for 
the singularity itself. In particular each side (1.5.2.2) defines a holomorphic
2-form outside of the origin for any $(i,j),i \ne j$ or $(k,l),k \ne l$.
In fact this form is just the residue of the log-form 
${{dz_1 \wedge ... \wedge dz_n} \over {F_1 \cdot \cdot \cdot F_{n-2}}}$
at non singular points (i.e. outside of the origin) of (1.5.2.1).
\par The adjoint ideal ${\cal A}_O $ in the local ring ${\cal O}_O$
of the origin  of a germ of complete intersection (1.5.2.1), according
to the description of the trace map 1.5.1 can be made explicit as 
 follows. Let
 $f: \tilde {\bf C}^n \rightarrow {\bf C}^n$ be an embedded resolution
 of (1.5.1.1). Then ${\cal A}_O$ consists of $\phi \in {\cal O}_O$ such
that $f^* (\phi \cdot  {{dw_i \wedge dw_j} \over {{\partial
(F_1,...,F_{n-2})} \over {\partial (w_1,...\hat w_i,..,\hat
w_j,..,w_n)}}})$ admits a holomorphic extension from $f^{-1}({\bf C}^n
-O)$ to   $\tilde {\bf C}^n$.
\par Similarly, the elements of 
$H^0({\cal A}(d_1+...+d_{n-2}-n-1) \subset H^0(\Omega_{{\bf
P}^n} (d_1+...+d_{n-2}))$ can be viewed as meromorphic forms with 
log singularities near non singular points of (1.5.1.1) having as
residue a 2-form on a non singular locus of $F$ and admitting a
holomorphic extension on $\tilde F$. 
 
 \bigskip 
\section{Ideals and polytopes of quasiadjunction.}
\subsection{Ideals of quasiadjunction.}
  Let $f$ be a germ of a {\it reduced}
 algebraic curve having
a singularity with $r$ irreducible branches at the origin of ${\bf C}^2$ 
near which it is given by local equation $f=f_1 (x,y) 
\cdot \cdot \cdot f_r(x,y)=0$. 
Let $\cal O$ be the local ring of the origin and $A$ be an ideal in
$\cal O$. 
\begin{dfn}
An ideal $A$ is called an ideal of quasiadjunction 
of $f$ with parameters $(j_1,..,j_r \vert m_1,...,m_r)$
($j_i,m_i$ are integers)
 if $A=\{\phi \in {\cal O} \vert z_1^{j_1}...z_r^{j_r} \phi
\in Adj V_{(m_1,f_1),...,(m_r,f_r)} \}$ where $V_{(m_1,f_1),,,(m_r,f_r)}$ 
is a germ at the origin of the complete intersection
 in ${\bf C}^{r+2}$ given by  the equations: 
  $$z_1^{m_1}=f_1(x,y),....,z_r^{m_r}=f_r(x,y).  \eqno (2.1.1)$$
An ideal of quasiadjunction is an ideal in $\cal O$ which is 
an ideal of quasiadjunction for some system of parameters.
\end{dfn}
\subsection{Basic Ideal.}
 Let ${A}(f_1,...,f_r) \subset {\cal O}$ be the ideal
generated by $${{(f_i)_x} \over {f_i}} f_1 \cdot f_2 \cdot \cdot \cdot
f_r, {{(f_i)_y} \over {f_i}} f_1 \cdot f_2 \cdot \cdot \cdot
f_r, (i=1,..,r), {{Jac({{(f_i,f_j)} \over {(x,y)}}))} \over {f_if_j}}
 f_1 \cdot \cdot \cdot f_r, (i,j=1,...,r,i \ne j) \eqno (2.2.1)$$ 
(we shall call it the basic ideal). 
\par Equating all polynomials (2.2.1) to zero yields a system of equations 
 having $(0,0)$ as
the only solution. Therefore ${\cal O}/A(f_1,...,f_r)$ is 
an Artinian algebra. 
\par Moreover for any set of parameters 
$(i_1,...,i_r \vert m_1,...,m_r)$ the corresponding ideal of 
quasiadjuction contains $A(f_1,...,f_r)$. Indeed, if
 $F_i=z_i^{m_i}-f_i(x,y)$, then up to sign
${{\partial (F_1,...,F_r)} \over {\partial (z_1,...,\hat z_i,..,z_r,x)}}=
z_1^{m_1-1} \cdot \cdot \cdot \hat{z_i^{m_i-1}} \cdot \cdot z_r^{m_r-1}
\cdot (f_i)_x$ and hence 
 $${{(f_i)_xf_1 \cdot \cdot 
\cdot \hat f_i \cdot \cdot \cdot f_r dz_i \wedge dy} \over 
 {{\partial (F_1,...F_r)} \over {\partial (z_1,..\hat z_i,...,z_r,x)}}}=
 {{(f_i)_xf_1\cdot \cdot \cdot \hat f_i \cdot \cdot \cdot 
 f_r dz_i \wedge dx} \over {z_1^{m_1-1} 
\cdot \cdot \cdot \hat z_i^{m_i-1} \cdot \cdot \cdot z_r^{m_r-1}(f_i)_x}}
=z_1 \cdot \cdot \cdot \hat z_i \cdot \cdot \cdot z_rdz_i \wedge dy
\eqno (2.2.2)$$
which is holomorphic on ${\bf C}^{r+2}$. Similarly, one sees that the 
2-forms corresponding to other generators of $A(f_1,..,f_r)$ coincide 
on $F_1=...=F_r=0$ with the forms admitting a holomorphic extension to 
${\bf C}^{r+2}$.  
\par In particular, there are only finitely many ideals of quasiadjunction.
\subsection{Ideals of quasiadjunction and polytopes.}
 Let ${\cal U}=\{ (x_1,..,x_r) \in {\bf R}^r \vert 
 0 \le x_i < 1 \}$ be the unit cube with coordinates  
corresponding to the components of a curve $\cal C$. Sometimes we shall 
denote this cube as ${\cal U} ({\cal C})$. If ${\cal C}'$ is formed by components
of $\cal C$ then we shall view ${\cal U}({\cal C}')$ as the face of 
 ${\cal U}({\cal C})$ given by $x_j=0$ where $j$ runs through indices
corresponding to components of $\cal C$ {\it not} belonging to ${\cal C}'$. 
\par By a {\it convex polytope} mean a subset of ${\bf R}^n$ 
\footnote{{only subsets of $\cal U$ will occur below}.} 
which is the convex hull of a finite set of points 
{\it with some faces possibly deleted}. By a {\it polytope} 
we mean a finite union of convex polytopes.
Class of polytopes in this sense is closed under finite 
unions and intersections.
A complement to a polytope within an ambient polytope is a
polytope. By {\it face of maximal dimension} of a polytope we mean the 
intersection of the polytope's boundary with the hyperplane
for which this intersection has the dimension equal to the dimension
of the boundary. A {\it face} of a polytope is an intersection of faces 
of maximal dimension. This is again a closed polytope.
\par Each $(x_1,..,x_r) \in {\cal U}$ defines the character 
$\chi (x_1,,,,,x_r)$ of
$\pi_1({\bf C}^2-{\cal C})$ such that 
$\chi (x_1,...,x_r)(\gamma_i)=exp(2 \pi \sqrt {-1} x_i)$. We shall call 
it the {\it exponential map}. We also  put $\allowbreak \overline {(x_1,...,x_r)}=
(1-x_1,...,1-x_r)$ and call this conjugation since 
$\chi (1-x_1,...,1-x_r)=\overline {\chi (x_1,...,x_r)}$. This map is an 
involution of the interior of 
$\cal U$: ${\cal U}^{\circ}=\{(x_1,...,x_r) \in {\cal U} \vert
x_i \ne 0, i=1,...r \}$.
\par Different arrays $(i_1,...,i_r\vert m_1,...,m_r)$ 
may define the same ideals of quasiadjunction. The next proposition 
describes when this is the case.
\begin{prop} Let ${A}$ be an ideal of quasiadjunction. 
Then there is a polytope $\bar {\Delta} (A)$, which is a open subset in
${\cal U}$, with the following property: for 
$(m_1,...,m_r) \in {\bf Z}^r$ and $(j_1,...,j_r) \in {\bf Z}^r,
0 \le j_i <m_i, i=1,...,r$  a holomorphic function 
$z_1^{j_1} \cdot \cdot \cdot z_r^{j_r} \phi$
belongs to the adjoint ideal of the germ of an abelian branched cover of 
the type $(m_1,..,m_r)$ of a neighborhood of the origin in ${\bf C}^2$ 
for any $\phi \in A$ if and only if  
 $({{j_1+1} \over {m_1}},...,{{j_r+1} \over {m_r}})
 \in \bar {\Delta}(A)$.
\end{prop}
{\bf Proof.}  Let $\rho: Y_f \rightarrow {\bf C}^2$
be an embedded resolution for the  singularity  
$f=f_1(x,y)\cdot \cdot \cdot f_r(x,y)=0$ at the origin. The 
complete intersection $V_{(m_1,f_1),...,(m_r,f_r)}$ (cf. (2.1.1)) 
provides a model with an 
isolated singularity of a branched abelian cover of a
neighborhood of the origin in ${\bf C}^2$ with $f=0$ as its branching locus.
Let $\pi_{(m_1,f_1),...,(m_r,f_r)}: 
V_{m_1,...,m_r} \rightarrow {\bf C}^2$ be the canonical projection. 
If $\bar V_{(m_1,f_1),...,(m_r,f_r)}$ is the normalization of 
$Y_{f} \times_{{\bf C}^2} V_{(m_1,f_1)...,(m_r,f_r)}$ then 
the projection 
 $\rho_{(m_1,f_1)...,(m_r,f_r)}: 
\bar V_{(m_1,f_1),...,(m_r,f_r)} \rightarrow V_{(m_1,f_1)...,(m_r,f_r)}$
on the second factor 
is a resolution of the singularity at the origin in the category of 
$V$-manifolds. We have the diagram: 
$$ \matrix{{ \atop  
 {\bar V_{(m_1,f_1),...,(m_r,f_r)}}}&{ {{\bar \pi_{(m_1,f_1),...,(m_r,f_r)}}
  \atop \longrightarrow} } 
& { \atop {Y_{f}}}\cr
{\Big \downarrow  \rho_{(m_1,f_1),...,(m_r,f_r)}} & & 
 {\Big \downarrow {\rho}} \cr
{\atop  V_{(m_1,f_1),...,(m_r,f_r)}} & {{\pi_{(m_1,f_1),...,(m_r,f_r)}} 
\atop {\longrightarrow}}& {\atop {\bf C}^2}
 \cr}  \eqno (2.3.1)$$
Let $E=\cup_k E_k$ be the exceptional locus of $\rho$. The
exceptional locus of $\rho_{(m_1,f_1),...,,(m_r,f_r)}$ is $\cup E_{k,l}$ 
where $E_{k,l}$ is a cover of $E_k$.
\par Let  $a_{k,i}=mult_{E_k} \rho^*f_i(x,y),(i=1,...,r)$, 
$c_k=mult_{E_k} \rho^* (dx \wedge dy)$. For $\phi \in {\cal O}$ 
we put  $f_k(\phi)=mult_{E_k} \rho^* (\phi)$. Finally let $g_{k,i}=g.c.d
(m_i,a_{k,i}),(i=1,...,r)$ and $s_k=
g.c.d. ({m_1 \over g_{k,1}},...{m_r \over
g_{k,r}}), (i=1,...,r)$. The form ${dx \wedge dy} \over {z_1^{m_1-1}
\cdot \cdot \cdot  z_r^{m_r-1}}$ is a non vanishing
form on $$V_{(m_1,f_1),..,(m_r,f_r)}-Sing V_{(m_1,f_1),..,(m_r,f_r)}$$
 (cf. (1.5.2.2)). If
${A}$ is an ideal of quasi-adjunction with parameters 
$(\bar j_1,...,\bar j_r \vert \allowbreak m_1,...,m_r)$ then the condition
 $\phi \in {A}$ is equivalent to the existence of a holomorphic extension of
 $\allowbreak \rho^*_{(m_1,f_1),...,(m_r,f_r)}({{z_1^{\bar j_1}...z_r^{\bar j_r}{\phi} dx \wedge dy} 
\over {z_1^{m_1-1}...z_r^{m_r-1}}})$ over the exceptional locus
$\cup E_{k,l}$ in a neighborhood of each point of $\cup E_{k,l}$
not belonging to $E_{k,l} \cap 
 E_{\bar k, \bar l}$ for any $(k,l),(\bar k,\bar l)$. This, in turn, is 
equivalent to:  
$$\Sigma_{i=1}^{i=r}(\bar j_i-m_i+1)mult_{E_{k,l}} 
\rho_{(m_1,f_1),....,(m_r,f_r)}^*(z_j)+ mult_{E_{k,l}}
 \rho_{(m_1,f_1),...,(m_r,f_r)}^* \phi+$$ $$+mult_{E_{k,l}} 
\rho_{(m_1,f_1),...,(m_r,f_r)}^* (dx \wedge dy) ) \ge 0 \eqno (2.3.2)$$
 for any pair of indices $(k,l)$.
On the other hand, we have the following equalities: 
$$mult_{E_{k,l}} \rho_{(m_1,f_1),...,(m_r,f_r)}^*(z_i)=
{{m_1 \cdot \cdot \cdot \hat m_i \cdot \cdot \cdot  m_r  \cdot
 a_{k,i}} \over  {g_{k,1} \cdot \cdot \cdot g_{k,r} \cdot s_k}},
mult_{E_{k,l}} \rho_{(m_1,f_1),...,(m_r,f_r)}^*(\phi)=$$ $$=
{{f_k(\phi) \cdot m_1 \cdot \cdot \cdot m_r}  \over {g_{k,1} \cdot \cdot \cdot
 g_{k,r} \cdot s_k }},$$
$$mult_{E_{k,l}} \rho_{(m_1,f_1),...,(m_r,f_r)}^*(dx \wedge dy)=
{{c_k \cdot m_1 \cdot \cdot \cdot m_r}  \over {g_{k,1} \cdot \cdot \cdot
 g_{k,r} \cdot s_k}}+
 {{ m_1 \cdot \cdot \cdot m_r}  \over {g_{k,1} \cdot \cdot \cdot
 g_{k,r} \cdot s_k }}-1 \eqno (2.3.3)$$
To see (2.3.3), we can 
 select  local coordinates $(u,v)$ on $Y_{f}$ near a point
belonging to a single component  $E_k$ in which the latter is given by
the  equation $u=0$. Then 
$\rho^*(f_i(x,y))=u^{a_{k,i}}\cdot \epsilon_i(u,v)$
where $\epsilon_i(u,v) (i=1,...,r)$ are units in the corresponding local ring.
The fiber product $Y_f \times_{{\bf C}^2} V_{(m_1,f_1),..,(m_r,f_r)} $
 is a subvariety in ${\bf C}^{r+2} \times_{{\bf C}^2} Y_f$
given by the equations: $$z_1^{m_1}=u^{a_{k,1}} \epsilon_{1}(u,v),...,
 z_r^{m_r}=u^{a_{k,r}} \epsilon_{r}(u,v)  \eqno (2.3.4)$$  
Each branch of (2.3.4) has the following local parameterization:
 $$u=t^{{m_1 \cdot \cdot \cdot m_r}  \over {g_{k,1} \cdot \cdot \cdot
 g_{k,r} \cdot s_k }},
 z_i=t^{{m_1 \cdot \cdot \cdot \hat m_i \cdot \cdot \cdot  m_r  \cdot
 a_{k,i}} \over  {g_{k,1} \cdot \cdot \cdot g_{k,r} \cdot s_k}}, 
i=1,...,r \eqno (2.3.5)$$
(exponents are chosen so that their greatest common
divisor will be equal to 1  and so that they will satisfy (2.3.4)).
This yields the first equality in (2.3.3). We have $$mult_{E_{k,l}} 
\rho^*_{(m_1,f_1),...,(m_r,f_r)}(\phi)=  
mult_{E_{k,l}} \bar \pi_{(m_1,f_1),...
,(m_r,f_r)} u^{f_k(\phi)}.$$ Hence the second equality in (2.3.3) follows 
from (2.3.5). 
\par Finally, since 
the map $\allowbreak  \bar \pi_{(m_1,f_1),...,(m_r,f_r)}$ 
is given locally by $(t,v) \rightarrow (t^{{m_1 \cdot \cdot \cdot m_r}
\over {g_{k,1} \cdot \cdot \cdot g_{k,r}}},v)$,
we have: $$\rho_{(m_1,f_1),...,(m_r,f_r)}^*(dx \wedge dy \vert 
V_{(m_1,f_1),...,(m_r,f_r)})=$$ 
$$=\rho^*_{(m_1,f_1),...,(m_r,f_r)} \circ 
\pi_{(m_1,f_1),...,(m_r,f_r)}^* (dx \wedge dy \vert {\bf C}^2)=$$
$$=\bar \pi^*_{(m_1,f_1),...,(m_r,f_r)}(u^{c_k}du \wedge dv)=
t^{{{c_k \cdot m_1 \cdot \cdot \cdot m_r} \over {g_{k,1} \cdot
\cdot \cdot g_{k,1} \cdot s_k}}+ {{m_1 \cdot \cdot \cdot m_r}
 \over {g_{k,1} \cdot \cdot \cdot g_{k,r} \cdot s_k}}-1}dt \wedge dv$$
which implies the last equality in (2.3.3).
\par Now it follows from (2.3.2) and (2.3.5) that
 $\phi \in {A}(\bar j_1,...,\bar j_r
 \vert m_1,...,m_r)$ if and only if for any $k$  the multiplicity $f_k(\phi)$ satisfies:
$$\Sigma_{i=1}^{i=r}
 (\bar j_i-m_i+1) {{m_1 \cdot \cdot \cdot \hat m_i \cdot \cdot \cdot m_r
 \cdot a_{k,i}} \over {g_{k,1} \cdot \cdot \cdot g_{k,r}s_k}}+
{{m_1 \cdot \cdot \cdot m_r \cdot f_k(\phi)} \over {g_{k,1} \cdot
\cdot \cdot g_{k,r} \cdot s_k}}+$$
$$+{{c_k \cdot m_1 \cdot \cdot \cdot m_r} \over {g_{k,1} \cdot \cdot \cdot
 g_{k,r} \cdot s_k}}+ {{m_1 \cdot \cdot \cdot m_r} \over {g_{k,1} \cdot
\cdot \cdot g_{k,r} \cdot s_k}}-1 \ge 0 \eqno (2.3.6)$$
For given $k$ let $f_k (A)$ be the minimal {\it integer} 
solution, with $f_k(\phi)$ considered as unknown, 
for this inequality and $\phi_k$ be such that 
$f_k(\phi_k)=f_k(A)$. In other words $\phi \in A$ if and only if  
 $f_k(\phi) \ge f_k(A)$.
 We have $f_k({A})=
\{\Sigma_{i=1}^{i=k} (a_{k,i}-(\bar j_i+1)
{{a_{k,i}}\over m_i}) -c_k-1 \}=
[(\Sigma_{i=1}^{i=k} (a_{k,i}-{{\bar j_i+1} \over m_i}a_{k,i})-c_{k} ]$ 
where  $\{ r \}$ (resp. $[r]$) denotes the smallest integer which 
is strictly greater than (resp. the integer part of) $r$.
We shall call $f_k(A)$ the
multiplicity of $A$ along $E_k$. This is the minimum of multiplicities along
$E_k$ of pull backs on $Y_{f}$ of elements of $A$.   
\par  The same calculation shows that $z_1^{j_1} \cdot \cdot \cdot z_r^{j_r}
 \phi$, where $\phi$
belongs to an ideal of quasiadjunction $A$, is in the adjoint ideal of 
$z_1^{m_1}=f_1(x,y),...,z_r^{m_r}=f_r(x,y)$ if and only if:
$$\Sigma_{i=1}^{i=r}(j_i-m_i+1){{m_1 \cdot \cdot \cdot m_r a_k} 
\over {g_{k,1} \cdot \cdot \cdot g_{k,r} \cdot m_i \cdot s_k}} 
 +{{m_1 \cdot \cdot \cdot m_r f_k(A)} \over {g_{k,1} \cdot \cdot \cdot 
 g_{k,r} \cdot s_k}}+$$ $$+{{c_k \cdot m_1 \cdot \cdot \cdot m_r} \over 
{g_{k,1} \cdot \cdot \cdot g_{k,r} \cdot s_k}}+{{m_1 \cdot \cdot \cdot m_r}
  \over {g_{k,1} \cdot \cdot \cdot g_{k,r} s_k}}-1 \ge 0 \eqno (2.3.7)$$
or equivalently:
$$\Sigma_{i=1}^{i=r} {{j_i+1} \over {m_i}}
a_{k,i} >\Sigma_{i=1}^{i=r} a_{k,i}-f_k(A)-c_k-1 \eqno (2.3.8)$$
Indeed, if (2.3.7) holds, then, since $f_k(\phi) \ge f_k(A)$ for 
any $\phi \in A$, one can replace $f_k(A)$ in (2.3.7) by $f_k(\phi)$.
This converts (2.3.7) into a necessary and sufficient condition
for $z_1^{j_1} \cdot \cdot \cdot z_r^{j_r} \phi$ to belong to 
the adjoint ideal of (2.1.1) (cf. derivation of (2.3.6)). 
Vice versa, for $\phi_k$ satisfying $f_k(\phi_k)=f_k(A)$,
the condition that
 $z_1^{j_1} \cdot \cdot \cdot z_r^{j_r} \phi_k $
is in the adjoint ideal of (2.1.1), is nothing else but (2.3.7).
\par The polytope $\bar {\Delta} (A)$ satisfying the conditions of the 
proposition is  the set of solutions of the 
inequalities:
 $$\Sigma_{i=1}^{i=r} x_i a_{k,i} >
\Sigma_{i=1}^{i=r} a_{k,i}-f_k(A)-c_k-1 \eqno (2.3.9)$$
 for all $k$. 
\smallskip \par \noindent
\subsubsection{Remarks.} 1. If $A_1$ and $A_2$ 
are ideals of quasiadjunction and $A_1 \subset A_2$ then 
$\bar \Delta (A_2) \subset \bar \Delta (A_1)$
\par \noindent
2.The polytope corresponding the basic ideal $A(f_1,..,f_r)$ (cf. 2.2)
is the whole unit cube $\cal U$.
\par \noindent 
3. ${\cal O}$ is considered as improper ``ideal'' of quasiadjunction 
since $A(m_1-1,...,m_r-1 \vert m_1,..,m_r)={\cal O}$.
\subsection{Local polytopes of quasiadjunction 
and faces of quasiadjunction.}
\smallskip \par \noindent
\subsubsection{}
\begin{dfn} 
We say that two points in the unit cube 
$\cal U$ are equivalent if the collections of polytopes $\bar \Delta (A)$
containing each of the points coincide. 
A ({\it local}) polytope of quasiadjunction $\Delta$ is an equivalence 
class of points with this equivalence relation. 
\end{dfn}
\begin{dfn}
A face of quasiadjunction is 
an intersection of a face (cf. (2.3)) of a local 
polytope of quasiadjunction 
and a (different) polytope of quasiadjunction. 
In particular, each face 
of quasiadjunction belongs
to a unique polytope of quasiadjunction.
\end{dfn}
\par For each face let us consider the system of equation defining 
the affine subspace of ${\bf Q}^r$ spanned by this face. One can 
normalize the system so that all coefficients of variables
and the free term are {\it integers} and the    
g.c.d. of non zero minors of maximal order is equal to 1. 
\smallskip \par \noindent 
\begin{dfn}The {\it order} of a face of 
quasiadjunction is the g.c.d of minors of maximal order in the 
{\it matrix of coefficients} in the normalized system of 
linear equations defining this face. 
\end{dfn}
This is the order of the torsion of the quotient of ${\bf Z}^{r}$
by the subgroup generated by the vectors having as coordinates 
the coefficients of variables in the equations of the face. 
In particular this integer is independent of the chosen 
normalized system of 
equations and depends only on the face of quasiadjunction. 
\par {\it The (local) ideal
of quasiadjunction corresponding to a face of quasiadjunction} 
is the ideal of 
quasiadjunction corresponding to the polytope of quasiadjunction
containing this face. 
\smallskip The ideal corresponding to a face of quasiadjunction has 
the form 
$A(j_1,...,j_r \vert m_1,...,m_r)$ where $({{j_1+1} \over {m_1}},
...,{{j_r+1} \over {m_r}})$ belongs to this face. 
\subsection{Examples.}  1. In the case of the branch with one component
the ideals of quasiadjunction correspond to the constants
of quasiadjunction (cf. [L2]). Recall that for a germ $\phi$ the rational
number $\kappa_{\phi}$ is characterized by the property 
$min \{ i \vert z^i \phi \in Adj(z^n=f(x,y)) \}=[\kappa_{\phi} \cdot n]$.
The ideal of quasiadjunction corresponding 
to $\kappa$ consists of $\phi$ such that $\kappa_{\phi} >\kappa$.
For example for the cusp $x^2+y^3$ the only non zero constant of 
quasiadjunction is $1 \over 6$, there are two polytopes of quasiadjunction
i.e. $\Delta'=\{ x \in [0,1] \vert 1 > x > {1 \over 6} \}$ and 
$\Delta''=\{ x \in [0,1] \vert { 1 \over 6}  \ge x \ge 0 \}$. 
$x= {1 \over 6}$ is the face of quasiadjunction and the corresponding ideal 
of quasiadjunction is the maximal ideal. 
For an arbitrary unibranched singularity 
we have $\bar \Delta(A_{\kappa})=\{x \vert x > \kappa \}$.
The order of the face $x=\kappa$ is equal to the order of the root 
of unity $exp(2 \pi i \kappa)$.
\par 2. Let us consider a tacnode, locally given by the equation:
 $y(y-x^2)$ (i.e. $f(x,y)=y, g(x,y)=y-x^2)$. 
Then the basic ideal is the
 maximal ideal and hence there is only one ideal of quasiadjunction. 
If $\cal M$ is the maximal ideal, then $\bar \Delta({\cal M})$ is the 
whole unit square. To determine 
the polytope $\bar \Delta({\cal O})$, note that after two blowups
we obtain an embedded resolution which in one of the charts looks like:
 $x=uv,y=u^2v$ where $u=0$ and $v=0$ are the exceptional curves. 
Hence for the component $u=0$ we obtain $a=b=2, f(\phi=1)=0,c=2$
and the corresponding polytope is $x+y > {1 \over2}$. 
The face of quasiadjunction is $x+y={1 \over 2}$ and the corresponding 
ideal is the maximal one. 
\par 3. For the ordinary singularity of 
multiplicity $m$: $(\alpha_1 x+\beta_1 y) \cdot \cdot \cdot 
 (\alpha_m x+\beta_m y)=0$ the basic ideal is ${\cal M}^{m-2}$ where
 ${\cal M} \subset {\cal O}$ is the maximal ideal of the local ring
at the origin. Since the resolution 
can be obtained by a single blow up, we have $a_{1,i}=1, c_1=1, i=1,...,m$ 
i.e. the polytope $\bar \Delta(A)$  of an ideal of quasiadjunction $A$ is:
  $$x_1+...+x_m > m-2-f_1(A) \eqno (2.5.1)$$ 
Since $f_1(\phi) \ge f_1(A)$ is equivalent to $\phi \in {\cal M}^{f_1(A)}$
i.e. the latter is the ideal corresponding to the polytope (2.5.1). 
The faces of quasiadjunction are $x_1+...+x_m=m-2-f_1(A)$ 
($f_1(A)=0,...,m-3$)
and the corresponding ideal of quasiadjunction is ${\cal M}^{f_1(A)+1}$.
\par Additional examples are discussed in [L7].
\subsection{Global polytopes and sheaves of quasiadjunction.} 
\smallskip \par 
Let ${\bf R}^r$, as in 2.3,  be  vector space coordinates of which 
are in one to one correspondence with the components of the curve 
${\cal C}=\cup_{i=1}^{i=r} C_i$. For a singular point $p$ of $\cal C$,
let ${\cal C}_p$ be the collection of components of $C$ passing through $p$.
Each polytope of quasiadjunction 
$\Delta_p \subset {\cal U}({\cal C}_p)$ of $p$ 
defines 
the polytope in ${\cal U}=\{(x_1,..,x_r) \vert 0 \le x_i \le 1 \} \subset 
{\bf R}^r$ consisting of points $\{ (x_1,...,x_r) \vert (x_1,..,x_r) \in
{\cal U},  (x_{i_1},...,x_{i_{r(p)}}) \in \Delta_p \}$ where 
$(i_1,...,i_{r(p)})$ are the coordinates corresponding to the
components of ${\cal C}_p$ (i.e. passing through the singularity $p$). 
We shall 
use the notation $\Delta_p({\cal U})$ for this polytope in $\cal U$. 
\subsubsection{Definition}
A type of a point in $\cal U$ is the 
collection of polytopes $\Delta_p({\cal U}) \subset {\cal U}$ 
to which this point belongs, where $p$ runs through all 
$p \in Sing {\cal C}$.
\par We call two points in $\cal U$ equivalent if they have the same
type. {\it A (global) polytope of quasiadjunction} is an equivalence
class of this equivalence relation.
\smallskip Global polytopes of quasiadjunction form a partition 
which is a refinement of every partitions of $\cal U$ defined by 
polytopes $\Delta_p({\cal U})$ corresponding to local 
polytopes of quasiadjunction of singularities of $\cal C$.
\smallskip \par \noindent
\subsubsection{Definition.}
 We shall call a face $\delta$ of  quasiadjunction
{\it contributing} if it belongs to a hyperplane $d_1x_1+...+d_rx_r=l$
where $d_1,...d_r$ are the degrees of the components of $\cal C$
corresponding to respective coordinates $x_1,...,x_r$. This hyperplane 
is called {\it contributing}, the integer $l=l(\delta)$ 
is called the {\it level}
of both the contributing hyperplane and contributing face.
The {\it order} of a global face of quasiadjunction is defined as 
in local case (cf. 2.4.2).
{\it A polytope of quasiadjunction is called contributing}
 if it contains
a contributing face. 
\par A point $(x_1,..,x_r) \in \delta$ is called {\it interior} 
$\C'$-point for a curve ${\cal C}'$ formed by components of $\cal C$ 
if $x_i \ne 1$ if and only if  $i$ corresponds to a component of ${\cal C}'$
and $(x_1,..,x_r)$ is in the interior of $\delta$.
\smallskip \par \noindent {\it Remarks}
\smallskip {\it  2.6.2.1} 
In the case when $r=1$, e.g. for an irreducible curve, an order of 
a global face of quasiadjunction is the order of a root of the 
Alexander polynomial. Indeed, for a constant of quasiadjunction $\kappa$,
$exp(2 \pi i \kappa)$ is a root of Alexander polynomial (cf. [LV]).
\par \noindent {\it 2.6.2.2.} The collection of orders of faces of 
quasiadjunction for reducible curves a priori cannot be determined 
just by the local types of all singular points. However it is 
{\it combinatorial} invariant of the curve in the sense that it 
depends only on local information about singularities and 
specification of components which contain specified singular points,
e.g. it is independent of the geometry of the set of singular points 
in ${\bf P}^2$. In the case of arrangements of lines this is combinatorial 
invariant in the common sense of the word.  
\smallskip \par \noindent 
\subsubsection{Definition.}  
The sheaf of ideals ${\cal A} (\delta) \subset {\cal O}_{{\bf P}^2}$ such that
$Supp ({\cal O}_{{\bf P}^2} /{\cal A}) \subset Sing {\cal C}$
is called {\it the sheaf of ideals of quasiadjunction} corresponding
to the face of quasiadjunction $\delta$ if the stalk ${\cal A}_p$ 
at each singular point  $p \in {\cal C}$ with local ring ${\cal O}_p$ 
is the ideal $A$ of quasiadjunction 
corresponding to the face $\Delta_p =\Delta \cap H_p$ with 
$H_p \subset {\bf R}^r$ being  given by $x_{i_j}=0$ where
$i_j$ are the coordinates corresponding to the components of $\cal C$
not passing through $p$.
\subsubsection{Examples.} 1.(cf. [L2]) For an irreducible curve of degree
$d$ with nodes and the ordinary cusps as the only singularities 
the global polytope of quasiadjunction coincide with the local one of the 
cusp. The only face of quasiadjunction is $x={1 \over 6}$. The contributing
hyperplane is given by $dx={d \over 6}$ and its level is $d \over 6$.
The sheaf of quasiadjunction corresponding to this face of quasiadjunction
is the ideal sheaf having stalks different from the local ring only 
at the points of ${\bf P}^2$ where the curve has cusps and the 
stalks at those points are the maximal ideals of the 
corresponding local rings.
\par \noindent 
2. Let us consider ${\cal C}$ which is an arrangement of lines.
For a point $P$ let $m_P$ denotes the multiplicity. 
We consider only points with $m_P>2$.
Each global face of quasiadjunciton is a solution of a system of equations:
 $$L_P: \ \ \ x_{i_1}+...+x_{i_m}=s_P \eqno (2.6.1)$$
where $s_P=1,...,m_P-2$ (cf. example 3 in 2.5). The indices
of variables $x$  
correspond to the lines of the arrangement and $x_i$ appears in 
$L_P$ if and only if it correspond to a line passing through $P$. 
Each system (2.6.1) corresponding to a face of quasiadjunciton
singles out a collection of vertices of the arrangement. 
This face is contributing if 
the equation
 $$x_1+...+x_r=k, \ \ \ (k \in {\bf N}) \eqno (2.6.2)$$
 is a linear combination of equations (2.6.1). 
The level of such contributing face is $k$.
Its order is the g.c.d of minors of maximal order in 
system (2.6.1).
\section{The first Betti number of an abelian cover.} 
\bigskip 
\par In this section we shall prove a formula for the irregularity
of an abelian covers of ${\bf P}^2$ branched over $\cal C$ in terms
of the polytopes of quasiadjunction introduced in the last section. More
precisely we shall calculate the multiplicity of a character of the Galois
group of the cover acting on the space $H^{1,0}(\tilde V_{m_1,..,m_r})$
of holomorphic 1-forms. We translate this into an information
about characteristic varieties of the fundamental group and consider
several examples of characteristic varieties for the fundamental groups 
of the complements to arrangements of lines.
 
\subsection{Statement of the theorem.}
 Let ${\cal C}=\cup_{i=1}^{i=r} C_i$ be a reduced 
curve $f(u,x,y)=f_1(u,x,y) \cdot \cdot \cdot f_r(u,x,y)$ 
with $r$ irreducible components and the degrees of
components equal to  $d_1,..,d_r$ and 
 $d=d_1+...+d_r$ be the total degree of $f(u,x,y)=0$.
Let $L_{\infty } $ be the line 
$u=0$ at infinity which, as above, we shall assume {\it transversal}
to $\cal C$ (cf.1.1). 
\par a) The irregularity of a desingularization  
$\widetilde V_{m_1,...,m_r}$ of an abelian cover of ${\bf P}^2$
branched over ${\cal C} \cup L_{\infty}$  
 and corresponding to the surjection  
$\pi_1 ({\bf P}^2-{\cal C} \cup L_{\infty}) \rightarrow  
{\bf Z}/m_1 {\bf Z} \oplus ... \oplus {\bf Z}/m_r{\bf Z}$ 
 is equal to 
 $$\Sigma_{{\cal C}'} (\Sigma_{\delta ({\cal C}')} N(\delta ({\cal C}'))
 \cdot  dim H^1 ({\cal A}_{\delta ({\cal C}')}(d-3-l(\delta({\cal C}))))) 
\eqno  (3.1.1)$$
where the summations are over all curves ${\cal C}'$ formed by the 
components of $\cal C$ and the contributing faces of 
quasiadjunction $\delta ({\cal C}')$ respectively. 
Here $l(\delta ({\cal C}'))$ is the level of the contributing 
face of $\delta({\cal C}')$ and $N(\delta({\cal C}'))$ 
is the number of interior ${\C}'$-points 
$({{i_1+1} \over m_1},...,{{i_r+1} \over m_r})$ in the contributing
face of $\delta({\cal C}')$.
\par b) Let $\chi_j$ be the character of 
${\bf Z}_{m_1} \oplus ...\oplus {\bf Z}_{m_r}$ taking on 
$(a_1,..,a_j,..,a_r)$ value 
\allowbreak  $exp(2 \pi {\sqrt {-1}}{{a_j} \over {m_j}})$. 
For a character $\chi$ of  ${\bf Z}_{m_1} \oplus ...\oplus {\bf
Z}_{m_r}$ let $$H^{1,0}_{\chi}(\widetilde V_{m_1,...,m_r})=\{ x \in 
H^{1,0}(\widetilde V_{m_1,...,m_r}), 
g \in {\bf Z}_{m_1} \oplus ...\oplus
 {\bf Z}_{m_r} \vert g \cdot x=\chi (g) \cdot x \}.$$
If $({{i_1+1} \over m_1},...,{{i_r+1} \over m_r})$ is an interior
$\cal C$-point (cf. 2.3.2) belonging to 
the contributing face $\delta$ then 
 $$dim H^1_{\chi_1^{i_1}...\chi_r^{i_r}}(\widetilde V_{m_1,...,m_r})=
  dim H^1 ({\cal A}_{\delta}(d-3-l(\delta))) \eqno (3.1.2)$$
\par c) Let $t_i =exp (2 \pi \sqrt {-1} x_i)$. For 
each contributing face $\delta$ belonging to ${\cal U}^{\circ}$ 
and its image $\bar \delta$ under the conjugation map (cf. 2.3), 
let $L_s(x_1,..,x_r)=\beta_s$ be the system of equations defining it
where $L_s (x_1,...x_r)$ is a linear form with integer coefficients
such that  g.c.d. of the minors of maximal order in the matrix
of coefficients is equal to $1$. 
Then the corresponding essential component of the
 characteristic variety of 
$\pi_1({\bf P}^2-{\cal C} \cup L_{\infty})$, which either has a positive
dimension or is a torsion point, 
is the intersection of cosets given by the equations: 
 $$exp (2 \pi \sqrt {-1} L_s)=exp(2\pi \sqrt {-1} \beta_s ) 
  \eqno (3.1.3)$$
written in terms of $t_i$'s. 
Vice versa, any essential component can be obtained in such way.
\bigskip \par Note that 
c) implies that the essential components of the characteristic 
varieties 
are Zariski's closures of the
images of the contributing faces under the exponential map.
Indeed, since g.c.d. of minors of coefficients in $L_s$ is $1$
the intersection of subgroups $exp (2 \pi \sqrt {-1} L_s)=1$ 
is connected and the closure of the image of the face of quasiadjunction
is Zariski dense in the translation of this connected component
given by (3.1.3).
\subsection{Proof of the Theorem}
 We shall start with the case when $m_i \ge d_i$ for $i=1,...,r$.
 Let ${\cal A} \subset {\cal O}_{{\bf P}^{r+2}}$ 
be the sheaf of adjoint ideals of the complete intersection
$V_{m_1,...,m_r} \subset {\bf P}^{r+2}$ given by the equations (cf. (1.3.1.1)):
 $$z_1^{m_1}=u^{m_1-d_1}f_1(u,x,y),...,
z_r^{m_r}=u^{m_r-d_r}f_r(u,x,y) \eqno (3.2.1)$$ 
 $V_{m_1,...,m_r}$ provides a model of an abelian branched cover 
of ${\bf P}^2$ branched over $f_1 \cdot \cdot \cdot f_r=0$ and the line
at infinity.  $V_{m_1,...,m_r}$ has 
isolated singularities at the points of (3.2.1) 
which are above the singularities of ${\cal C}$ in 
${\bf P}^2-L_{\infty}$. The action of the Galois group of the 
cover is induced from the action of the product of groups of roots of unity
$\mu_{m_1} \times \cdot \cdot \cdot \times \mu_{m_r}$ on the 
${\bf P}^{r+2}$ via multiplication of corresponding $z$-coordinates. 
\par Let $H$ be the set of common zeros of  
$z_1,...,z_r \in H^0({\bf P}^{r+2}, {\cal O}(1))$
and ${\cal A}_{i_1,...,i_r}$ be the subsheaf of ${\cal O}_{{\bf P}^{r+2}}$
 germs of section  product  of which with  $z_1^{i_1} \cdot \cdot \cdot
 z_r^{i_r}$ belongs to  ${\cal A}$. The action of $\mu_{m_1} \times 
\cdot \cdot \cdot \times \mu_{m_r}$ on ${\bf P}^{r+2}$ induces the action on 
${\cal A}_{i_1,...,i_r}$.
\par  
Let ${\cal J}_H$ be the ideal 
sheaf of the plane $H \subset {\bf P}^2$. 
We have the following $\mu_{m_1} \times
\cdot \cdot \cdot \times \mu_{m_r}$-equivariant sequence:
 $$0 \rightarrow {\cal A}_{i_1,...,i_r}((m_1-i_1)+...+(m_r-i_r)-r-3) 
\otimes {\cal J}_H 
\rightarrow {\cal A}_{i_1,...,i_r}((m_1-i_1)+...+(m_r-i_r)-r-3)$$
$$ \rightarrow {\cal A}_{i_1,..,i_r}((m_1-i_1)+...(m_r-i_r)-r-3) \vert _H 
\rightarrow 0 \eqno (3.2.2)$$
Let $$F(i_1,...,i_r)=dim H^1 ({\cal A}_{i_1,...,i_r}((m_1-i_1)+...+
(m_r-i_r)-r-3));  $$
 $$F_{\chi}(i_1,...,i_r)=dim \{ x \in H^1({\cal A}_{i_1,...,i_r}
 ((m_1-i_1)+...+(m_r-i_r)-r-3) \vert g \cdot x= \chi (g) x,$$
 $$\forall g \in \mu_{m_1} \times \cdot \cdot \cdot \mu_{m_r} \} \eqno (3.2.3)$$
In particular $F(0,...,0)$ is the irregularity of a nonsingular model
of $V_{m_1,...,m_r}$.  
\smallskip
\par  {\it Step 1. Degree of the curves in the linear system}
 $H^0({\cal A}_{i_1,...,i_r}((m_1-i_1)+...+(m_r-i_r)-r-3) \vert _H)$.
 Let us calculate the multiplicity of the line 
$L_{\infty} :z_1=...=z_r=u=0$ as the fixed component of the curves 
in the linear system cut on $H$ by the hypersurfaces in the linear system 
$ H^0({\cal A}_{i_1,...,i_r}((m_1-i_1)+...+(m_r-i_r)-r-3))$. 
This multiplicity 
is the smallest $k$ such that $u^k$ belongs to the latter system 
of hypersurfaces.
In appropriate coordinates $(z_1,...,z_r,u,v)$ at a  point $P$ 
of this line outside of $L_{\infty} \cap {\cal C}$ (i.e. we have 
 $f_1(P) \cdot \cdot \cdot f_r(P) \ne 0$) the local 
equation of $V_{m_1,...,m_r}$ is $z_1^{m_1}=u^{m_1-d_1}, ..., z_r^{m_r}
 =u^{m_r-d_r}$. 
Let $$l=l.c.m.({{m_1} \over {m_1-d_1}}(m_1-d_1) \cdot \cdot \cdot (m_r-d_r),
..., $$ $$...{{m_r} \over {m_r-d_r}} (m_1-d_1) \cdot 
\cdot \cdot (m_r-d_r),(m_1-d_1) \cdot \cdot \cdot (m_r-d_r))$$
 Then each branch of the 
 normalization of $V_{m_1,...,m_r}$ has the parameterization $(t,v)$
 such that:
 $$z_1=t^{l(m_1-d_1) \over {m_1 (m_1-d_1)
  \cdot \cdot \cdot (m_r-d_1)}},...,z_r=t^{l(m_r-d_r) \over {m_r(m_1-d_1)
 \cdot \cdot \cdot (m_r-d_r)}}, u=t^{l \over
 {(m_1-d_1) \cdot \cdot \cdot (m_r-d_r)}}$$
Therefore the pull back of the form ${{z_1^{i_1} \cdot \cdot \cdot 
z_r^{i_r} u^kdu \wedge dv} \over
 {z_1^{m_1-1} \cdot \cdot \cdot z_r^{m_r-1}}} \vert _{V_{m_1,...,m_r}}$
 to the $(t,v)$ chart is regular if and only if 
 $$ k>\Sigma_{j=1}^{j=r} (m_j-d_j-(i_j+1))+{{d_j(i_j+1)} \over {m_j}}
 -1 \eqno (3.2.4)$$
The smallest $k$ which satisfies this inequality, i.e. the multiplicity of 
the line $u=0$ as the component of a generic curve from  
$H^0({\cal A}_{i_1,...,i_r}((m_1-i_1)+...(m_r-i_r)-r-3) \vert _H)$, 
is equal to 
$$\Sigma_j (m_j-d_j-(i_j+1))+[\Sigma_j {{d_j(i_j+1)} \over m_j}]
 \eqno (3.2.5) $$
A consequence of this is that the degree of the moving curves in the 
linear system $\allowbreak 
H^0({\cal A}_{i_1,...,i_r} ((m_1-i_1)+...(m_r-i_r)-r-3) 
\vert H)$  is equal to $\Sigma_j d_j -3-[\Sigma_j 
{{d_j(i_j+1)} \over m_j}]$ and therefore the moving curves belong
to the linear system $H^0 ({\cal A}_{\Delta}((\Sigma_jd_j)-3-
[{{d_j(i_j+1)} \over {m_j}}])$ where $\Delta$ is the
polytope of quasiadjunction containing 
$({{i_1+1} \over m_1},...,{{i_r+1} \over m_r})$. 
In fact the moving curves form a complete system since the cone over
any curve in $H^0({\cal A}_{\Delta}(\Sigma_j d_j-3-[\Sigma_j {{d_j(i_j+1)}
 \over {m_j}}])$ belongs to $H^0({\cal A}_{i_1,...,i_r} (\Sigma_j
 (m_j-i_j)-r-3))$.
\smallskip 
\par {\it Step 2. A recurrence relation for $F(i_1,...,i_r)$ 
and $F_{\chi}(i_1,...,i_r)$.}
Let $s(i_1,...,i_r)=dim H^1({\cal A}_{\Delta}(\Sigma_jd_j-3-
[\Sigma_j{{d_j(i_j+1)} \over {m_j}}])$ where $\Delta$ is the polytope
of quasiadjunction of ${\cal C}$ containing $({{i_1+1} \over m_1},...,
{{i_r+1} \over m_r})$ and $\epsilon_{\chi}(i_1,...,i_r)=1$ 
(resp. $0$)  
if $\chi=\chi_1^{i_1-m_1+1} \cdot \cdot \cdot \chi_r^{i_r-m_r+1}$ 
(resp. otherwise).
We claim the following recurrence:
$$F(i_1,...,i_r)=s(i_1,...,i_r)+\Sigma_{l=1}^{l=r} 
(-1)^{l+1}\Sigma_{i_{j_1} <...<i_{j_l}}
 F(...,i_{j_1}+1,...,i_{j_l}+1,...); $$
$$F_{\chi}(i_1,...,i_r)=\epsilon_{\chi}(i_1,...,i_r)s(i_1,...,i_r)+$$
$$ \Sigma_{l=1}^{l=r}(-1)^{l+1} \Sigma_{i_{j_1} <...<i_{j_l}} 
F_{\chi(\chi_{j_1} \cdot \cdot \cdot \chi_{j_l})^{-1}}
(...,i_{j_1}+1,...,i_{j_l}+1,...) \eqno(3.2.6)$$ 
Equivalently the first of equalities (3.2.6) can be written as
$$s(i_1,...,i_r)=\Sigma_{l=0}^{l=r} 
(-1)^{l}\Sigma_{i_{j_1} <...<i_{j_l}}
 F(...,i_{j_1}+1,...,i_{j_l}+1,...)$$
and similarly for the second.
This identity will be derived from the following. For $h$ such that 
$1 \le h \le r$ let $$F(i_1,.,,i_r \vert q_1,...,q_h)=
dim H^1({\cal A}_{i_1,..,i_r}((m_1-i_1)+...(m_r-i_r)-r-3) \vert_{H_{q_1}
\cap ...H_{q_h}})$$ where $H_s$ is the hyperplane $z_s=0$ in ${\bf
P}^{r+2}$ while for $h=0$ we let 
$F(i_1,...,i_r \vert \emptyset)=F(i_1,...,i_r)$.
 In particular  $s(i_1,..,i_r)=F(i_1,...i_r \vert 1,...,r)$.
\par Similarly one defines $F_{\chi}(i_1,...,i_r \vert q_1,...,q_h)$.
We shall prove by induction over $h$:
 $$F(i_1,...,i_r \vert q_1,...,q_h)=$$
 $$\Sigma_{l=0}^{l=h} (-1)^l \Sigma_{i_{j_1}<...<i_{j_l},
 i_{j_1},...,i_{j_l} \subset (q_1,...,q_h)}
  F(...,i_{j_1},...,i_{j_l},...)$$
 $$F_{\chi}(i_1,...,i_r \vert q_1,...,q_{h})=$$ $$\Sigma_{l=0}^{l=h} 
(-1)^{l}
\Sigma_{i_{j_1} <...<i_{j_l},(i_{j_1},...,i_{j_l}) \subset (q_1,...,q_h)}
 F_{\chi \cdot \chi_{j_1}^{-1}
 \cdot \cdot \cdot \chi_{j_l}^{-1}}
 (...,i_{j_1}+1,...,i_{j_l}+1,... ) \eqno(3.2.7)$$
The identity (3.2.6) is a special case of (3.2.7) when $q_i=i$.
For any $(i_1,...,i_r \vert q_1,...,q_h),(h \ge 0)$ from the exact sequence
(in which the left map is the multiplication by $z_{q_{h+1}}$):
$$0 \rightarrow {\cal A}_{...,i_{q_{h+1}+1},...}
(...+(m_{q_{h+1}}-i_{q_{h+1}}-1)+...-r-3)\vert_{H_{q_1} \cap
... \cap H_{q_h}} \rightarrow $$
$${\cal A}_{i_1,..,i_{r}} ((m_1-i_1)+...+(m_r-i_r)-r-3) \vert_{H_{q_1} 
\cap ... \cap H_{q_h}} \rightarrow $$
$$ {\cal A}_{i_1,..,i_r}((m_1-i_1)+...+(m_r-i_r)-r-3)
\vert_{H_{q_1,...,q_h,q_{h+1}}} \rightarrow 0 \eqno (3.2.8)$$
we obtain $$F(i_1,...,i_r \vert q_1,...,q_h,q_{h+1})=-
F(i_1,...,i_{q_{h+1}}+1,...i_r
\vert q_1,..,q_h)+F(i_1,...,i_r \vert q_1,..,q_h) $$
$$F_{\chi}(i_1,...,i_r \vert q_1,...,q_h,q_{h+1})=
 -F_{{\chi} \cdot \chi_{q_{h+1}}^{-1}} 
 (i_1,...,i_{q_{h+1}+1},...,i_r \vert q_1,...,q_h)+$$
 $$ F_{\chi}(i_1,...,i_r 
 \vert q_1,..,q_h).  \eqno (3.2.9)$$ Indeed, the 
map $$H^0 ({\cal A}_{i_1,..,i_r}((m_1-i_1)+...+(m_r-i_r))-r-3)
\vert_{H_{q_1,...,q_h}} \rightarrow$$ $$\rightarrow 
H^0 ({\cal A}_{i_1,..,i_r}((m_1-i_1)+...+(m_r-i_r))-r-3)
\vert_{H_{q_1,...,q_h,q_{h+1}}} \eqno (3.2.10)$$
is surjective because  the cone in
$H_{q_1} \cap ... \cap H_{q_h}$ over the
hypersurface in $H^0({\cal A}_{i_1,..,i_r}((m_1-i_1)+...+(m_r-i_r)-r-3)
\vert_{H_{q_1,...,q_h,q_{h+1}}})$ belongs to
 $H^0({\cal A}_{i_1,..,i_r}((m_1-i_1)+...+(m_r-i_r)-r-3)
\vert_{H_{q_1,...,q_h}})$.
 Moreover for $q_i=i, i=1,...,r$ we have $F_{\chi}(i_1,...,i_r \vert
 1,...r)=\epsilon_{\chi}(i_1,...,i_r)s(i_1,...,i_r)$ since to 
$\phi(x,y) \in H^0({\cal A}_{i_1,...,i_r})$ corresponds the form 
$\psi=z_1^{m_1-i_1-1} \cdot \cdot \cdot z_{r}^{m_r-i_r-1} \pi^* \phi$
holomorphic on $\widetilde V_{m_1,...,m_r}$ and satisfying:
$g^*(\psi)=\chi_1^{i_1-m_1+1} \cdot \cdot \cdot \chi_r^{i_r-m_r+1}\psi$.
This shows that (3.2.7) is valid for $h=1$ and that validity of (3.2.7)
for the array $(q_1,..,q_{h+1})$ provided it is valid for all $(q_1,..,q_h)$.
\smallskip 
\par {\it Step 3. An explicit formula for} $F(i_1,...,i_r)$.
Let ${\cal C}(j_1,...,j_s)=C_{j_1} \cup ...\cup C_{j_s}$ be a
curve formed by a union of the components of ${\cal C}$ and 
let  
$$F_{{\cal C}(j_1,...,j_s)}(i_1,...,i_{s})= dim H^1({\cal A}
({\cal C}(j_1,...,j_s))_{i_1,...,i_s}((m_{j_1}-i_1)+....+
(m_{j_s}-i_s)+s-3).$$ Note that
 $$i_{j}=m_{j}-1 (j \ne j_1,...j_s) \Rightarrow  
F(i_1,...,i_r)=F_{{\cal C}(j_1,...,j_s)}(i_1,...,i_s) \eqno(3.2.11) $$ 
since the local conditions defining both sheaves
coincide (indeed: $${{z_1^{i_1} \cdot \cdot \cdot z_{j_k}^{m_{j_k}-1}
 \cdot \cdot 
 z_r^{i_r} dx \wedge dy} \over
{{z_1}^{m_1-1}...z_r^{m_r-1}}}={{z_1^{i_1}\cdot \cdot \cdot 
\hat {z_{i_k}^{m_{i_k}-1}}
\cdot \cdot \cdot z_r^{i_r}dx \wedge dy} \over {z_1^{m_1-1} \cdot \cdot 
\cdot \hat {z_{i_{j_k}}^{m_{j_k}-1}}
 \cdot \cdot \cdot z_r^{m_r-1}}})$$  
as well as the degrees of the curves in the corresponding linear 
systems. Moreover $$F_{{\cal C}(j_1,...,j_s)}(0,...,0)$$ is the irregularity 
of the cover of ${\bf P}^2$ branched over ${\cal C}(j_1,...,j_s)$
and having the ramification index $m_i$ over the component
$C_i (i=j_1,...,j_s)$. 
\par We solve the recurrence  
relation (3.2.6) subject to the  "initial condition"  (3.2.11).
It is convenient to view each relation (3.2.6) as the one connecting the 
values of the function defined at the vertices of the integer lattice 
in the parallelepiped $0 \le x_i \le m_i, (i=1,...,r)$. Each equation 
connects the values of this function at the vertices of a parallelepiped   
with sides equal to 1. It is clear that the sum of all equations (3.2.6)
yields:   
$$F(0,...,0)=\Sigma_{0 \le i_s < m_s-1}s(i_1,...,i_r) + 
\Sigma_{l=1}^{l=r-1} \Sigma_{(j_1<...<j_l)}
F_{{\cal C}(j_1,...,j_l)}(0,...0) $$
 $$F_{\chi}(0,...,0)=\Sigma_{0 \le i_s <m_s-1}\epsilon_{\chi}(i_1,..,,i_r)
 s(i_1,...,i_r)+
 \Sigma_{l=1}^{l=r-1} \Sigma_{(j_1<...<j_l)} F_{{\cal C}(j_1,...,j_l)_{\chi}}
(0,...,0) \eqno (3.2.12)$$
\par {\it Remark. Alternative derivation of (3.2.12).}
\smallskip \par  Sheaves 
${\cal A}_{\Delta}(d_1+...+d_r-r-2-[\Sigma_j {{d_j(i_j+1)} \over {m_j}}])$
admit the following interpretation also yielding (3.2.12). 
Let us consider the following global version of the diagram (2.3.1):
  $$\matrix {{\bar V_{m_1,...,m_r}} & {\bar \pi} \atop \rightarrow & 
     {Y_{\cal C}} \cr
 {\downarrow {\bar \rho }} &   & {\downarrow {\rho}} \cr
            {V_{m_1,...,m_r}} & \pi \atop \rightarrow  & {\bf P}^2 \cr} 
\eqno (3.2.13)$$
Here $\rho: Y_{\cal C} \rightarrow {\bf P}^2$ is an embedded resolution 
of singularities of ${\cal C}$ which are worse than nodes, 
$\bar V_{m_1,...,m_r}$ is the normalization of $V_{m_1,...,m_r} 
\times_{{\bf P}^2} Y_{\cal C}$ and ${\bar \pi}, {\bar \rho}$ are the 
obvious projections.    
Let $$\bar \pi_* ({\cal O}_{\bar V_{m_1,...,m_r}})=\oplus  
{\cal L}_{\chi_1^{i_1} \cdot \cdot \cdot \chi_r^{i_r}}^{-1} \eqno (3.2.14)$$
 be the decomposition by the characters
 of the Galois group acting on 
$\bar \pi_*({\cal O}_{\bar V_{m_1,...,m_r}})$.
Then we have:
$${\cal A}_{\Delta}(\Sigma_j d_j-3-[\Sigma_j {{d_j(i_j+1)} 
\over {m_j}}])=\rho_* ({\bar \pi}_* (\Omega_{\bar V_{m_1,...,m_r}}) 
\otimes {\cal L}_{\chi ({m_1-(i_1+1),...,m_r-(i_r+1)})}) \eqno (3.2.15)$$
where $\Delta$ is the polytope of quasiadjunction containing 
 $({{i_1+1} \over {m_1}},...,{{i_r+1} \over {m_r}})$.
Indeed it follows from (2.3.8) that a germ $\phi$ of a holomorphic function 
belongs to the sheaf in the left side of (3.2.15)
if and only if the order of $\phi$ along an exceptional curve $E_k \subset Y_{\cal C}$
satisfies: $ord_{E_k} \phi \ge 
\Sigma a_{k,j} ({{m_j-(i_j+1)} \over m_j}-c_k)$ and the sheaf on the left
is a subsheaf of ${\cal O}_{{\bf P}^2}([\Sigma_j d_j({{m_j-(i_j+1)} \over
m_j})])$ with the quotient having a zero-dimensional support.
One readily sees that the sheaf on the right has the same local description. 
This identity also implies (3.2.12) as follows from Serre's duality and
(3.2.14).
\smallskip 
\par {\it Step 4. A vanishing result.}
\par  If $\Delta$ is a polytope of quasiadjunction, 
 ${\Xi}_k=\{(x_1,...,x_r) \in {\cal U} \vert k \le d_1x_1+...+
d_rx_r < k+1 \} $ and $k$ is such that $\Delta \cap {\Xi}_k \ne
\emptyset$ then $$H^1({\cal A}_{\Delta}(d-r-2-k))=0 \eqno (3.2.16)$$
unless $\Delta$ is a contributing polytope of quasiadjunction
and $\Delta \cap \Xi_k$ is a face of quasiadjunction.
\par If $\Delta$ isn't contributing (cf. 2.6.2), then the intersection 
 $\Delta \cap \Xi_k$ has a positive volume. If $X(n)$ is the number of 
points $({i \over n},...,{i \over n})$ in the latter, it follows from
(3.2.12) that we have $b_1({\cal C},n) \ge 
dim H^1({\cal A}_{\Delta}(d-r-2-k)) \cdot X(n)$. We have $X(n) > C
 \cdot n^r$ for some non zero constant  $C$. Therefore we get 
contradiction with Corollary (1.3.3)
 unless $dim H^1({\cal A}_{\Delta}(d-r-2-k))=0$.
\par {\it Step 5. End of the proof.}
 Step 4 and the formula (3.2.12) give a) and b) of the theorem
 in the case $m_i \ge d_i$ for $i=1,..,i_r$. 
\par If $\chi$ is a character of $\allowbreak 
{\bf Z}_{m_1} \oplus ... \oplus {\bf Z}_{m_r}$
acting on $H^0 (\Omega^1_{\tilde V_{m_1,...,m_r}})=
H^{1,0}(\tilde V_{m_1,...,m_r})$  then $\bar \chi$ 
is a character with eigenspace of the same dimension for the action of
${\bf Z}_{m_1} \oplus ...\oplus {\bf Z}_{m_r}$ on $H^{0,1}(\tilde V_{m_1,..,m_r})$.
Hence  part b) and Sakuma formula (cf. 1.3.2) 
imply that a points $(.... {{i_j+1} \over m_j}...)$ belongs
to a contributing face of $\Delta$ or its conjugate 
if and only if  $\allowbreak (....,exp(2 \pi \sqrt{-1}
{{i_j+1} \over m_j}),...)$
belongs to $i$-th characteristic variety with 
$i=dim H^1({\cal A}_{\Delta}(d_1+...+d_r-r-3-k(\Delta))$.
Since a characteristic variety is a translated by a point 
of finite order subtorus (cf. 1.4.2) this implies c). 
Now the remaining cases of the formula a) follows from 
 Sakuma's result (1.3.2.2).
\subsection{Examples.} In 2.6.4 we did describe 
systems of equations
for faces of quasiadjunction in the case of arrangements of lines.
To determine if a set of solutions of the system corresponding to a face
$\delta$ actually corresponds to a component 
of characteristic variety one should 
\par a) calculate the superabundance
(3.1.2) of the corresponding linear system and 
\par b) decide the ``amount 
of translation'' i.e. to normalize the system of equations so that 
the g.c.d. of minors of the left hand sides of  (2.6.1) will be equal to
 one.
\par In any event, if superabundance is not zero, then clearly 
the component of characteristic variety will be a connected component of the 
subgroup given by the equations: $exp(L_P)=1$ with $P$ running 
through all vertices singled out by the face of quasiadjunction.
\smallskip \par \noindent {\it Example 1.} 
Let us calculate the irregularity of the abelian
cover of  ${\bf P}^2$ branched over the arrangement $L: uv(u-v)w=0$ and
 corresponding  to the homomorphism $H_1 ({\bf P}^2-L) =
 {\bf Z}^3 \rightarrow ({\bf Z}/n{\bf Z})^3$. The only nontrivial ideal of
 quasiadjunction is the maximal ideal of the local ring with
 corresponding polytope of quasiadjunction: $x+y+z >1$. Hence the
 irregularity of the abelian cover is $Card \{(i,j) \vert 0 <i <n, 0 <j
 <n, {i \over n}+{j \over n}+{k \over n}=1 \}
 \cdot dim H^1 ({\cal J}(3-3-1))$ where
${\cal J}=Ker {\cal O } \rightarrow {\cal O}_P$ where $P: u=v=0$.
$\cal J$ has the following Koszul resolution: $$0 \rightarrow {\cal O} (-2)
 \rightarrow
{\cal O}(-1) \oplus {\cal O} (-1) \rightarrow {\cal J} \rightarrow 0$$ which
 yields $H^1 ({\cal J}(-1))=H^2({\cal O}(-3))={\bf C}$. Now the counting
 points on $x+y+z=1$ yields ${n^2-3n+2} \over 2$ as the irregularity
of the abelian cover.
\smallskip \par \noindent
{\it Example 2.} Let us consider the arrangement formed by 
the sides of an equilateral triangle ($x_1,x_2,x_3$) and its medians
($x_4,x_5,x_6)$ arranged so that the vertices are the intersection
points of $(x_1,x_2,x_4),(x_2,x_3,x_5)$ and $x_3,x_1,x_6$ respectively
(Ceva arrangement cf. [BHH]).
It has 6 lines, 4 triple  and 3 double points. 
The polytopes of quasiadjunction are the connected components of the partition  
of ${\cal U}=\{(x_1,..,x_6) \vert 0 \le x_i \le 1, i=1,...6 \}$ by the
hyperplanes:
 $$x_1+x_2+x_4=1,x_2+x_3+x_5=1, x_3+x_1+x_6=1, x_4+x_5+x_6=1 \eqno (3.3.1)$$ 
The only face of a polytope of quasiadjunction which belongs to a 
hyperplane $H_k: x_1+x_2+x_3+x_4+x_5+x_6=k, k \in {\bf Z}$
is formed by set of solutions of the system of {\it all} 4 equations
(3.3.1). This face belongs to $H_2$ and is the only contributing face.
Hence the irregularity is equal
to $N \cdot dim H^1({\cal J} (6-3-2))$ where $N$ is the number of
solutions  (3.3.1) of the form $x_i={j \over n}$. To calculate
 $dim H^1({\cal J} (6-3-2))$ notice that 4 triple points form a complete 
intersection of two quadrics. 
This yields $H^1 ({\cal J} (1))=H^2 ({\cal O} (-3))={\bf C}$. 
\par It follows from (3.3.1) that the only essential torus is
a component of subgroup:
$$t_1t_2t_4=1, t_2t_3t_5=1,t_1t_3t_6=1, t_4t_5t_6=1 \eqno (3.3.2)$$
This subgroup has two connected components: 
$$(u,v,u^{-1}v^{-1},u^{-1}v^{-1},u,v), (-u,-v,-u^{-1}v^{-1},
u^{-1}v^{-1},u,v), \ \ \ u,v \in {\bf C}^* \eqno (3.3.3)$$
The second component is a translation of the first by $(1,1,1,-1,-1,-1)$,
a point of order 2. Since (3.3.1) admits an integral solution 
image under the exponential map  of the contributing face does contains 
trivial character
and hence the subgroup in (3.3.3) is the essential torus.
\par There are also 4 nonessential tori corresponding to each of 
triple points:
   $$t_1t_2t_3=1, t_i=1, i \ne 1,2,3 $$
   $$t_5t_2t_3=1, t_i=1, i \ne 5,2,3 $$
   $$t_4t_6t_3=1, t_i=1, i \ne 4,6,3 $$
   $$t_4t_5t_6=1, t_i=1, i \ne 4,5,6  \eqno (3.3.3) $$
Let us consider the abelian cover of ${\bf C}^2$ corresponding to the
homomorphism $H_1({\bf C}^2-C) \rightarrow 
({\bf Z}/n{\bf Z})^6/{\bf Z}/n{\bf Z}$ (embedding of the quotiented  
subgroup is diagonal). Then each of five tori contributes 
the same number into irregularity equal to ${n^2-3n+2} \over 2$,
i.e. the irregularity of the abelian cover of ${\bf CP}^2$ is 
$5{{n^2-3n+2} \over 2}$ (e.g., for $n=5$ the irregularity is $30$, cf. [I]).
\par \noindent
{\it Example 3.} Let us calculate the characteristic varieties of the union of 9 lines which are dual to nine inflection points
on a non singular cubic curve $C \subset {\bf P}^2({\bf C})$.
 This arrangement in ${{\bf P}^2}^*({\bf C})$ has 12 triple points
corresponding to 12 lines determined by the pairs of the inflection 
points of $C$. One can view inflection points of $C$ as the points
of ${\bf F}^2_3$ (${\bf F}_3$ is the field with 3 elements)
 i.e. as the points of the affine part in a projective plane
 ${\bf P}^2({\bf F}_3)$. The triple points of this arrangement then 
can be viewed as lines in ${\bf P}^2({\bf F}_3)$ different from the 
line at infinity (i.e. the complement to the chosen affine plane). 
In dual
picture one identifies triple points of this arrangement
with points of the dual plane ${{{\bf P}^2}^*}({\bf F}_3)$ 
different from a fixed point $P$ corresponding to the line at infinity.
Then the lines of this arrangement in ${{\bf P}^2}^*({\bf C})$ are 
identifies with the lines in ${{\bf P}^2}^*({\bf F}_3)$ 
not passing through the fixed point $P$. 
\par Each essential component corresponds to a collection of 
vertices $\cal S$ (cf. (2.6.4), example 2). 
The structure of the system of equations (2.6.1) 
shows that ${{\vert {\cal S} \vert} \over k}={r \over m}=3$.
Hence one has either: 
\par a) $\vert {\cal S} \vert=3,k=1$ or
\par b) $\vert {\cal S} \vert=6, k=2$ or
\par c) $\vert {\cal S} \vert=9, k=3$ or
\par d) $\vert {\cal S} \vert=12, k=4$. 
\par Cases a) and b) will not  define non empty tori since in this case
$r^2 > 9 {\vert {\cal S} \vert}$ (cf. corollary 4.1).
\par In the case c) each collection $\cal S$ is determined by one of 
4 choices of a line $\ell$ through $P$ and consists of 9 points in 
${{\bf P}^2}^*({\bf F}_3)$ in the complement to the chosen line. 
In this case 
the corresponding homogeneous system has rank 7 i.e. 
a 2-dimensional space of solutions. Moreover,    
$dim H^1 ({{\bf P}^2}^* ({\bf C}), {\cal I}(9-3-3))=1$ since the points
on  ${{{\bf P}^2}^*} ({\bf C})$ corresponding to 9 points in 
${{\bf P}^2}^*({\bf F}_3)$ in the complement to 
a line $p \subset  {{\bf P}^2}^*({\bf F}_3)$
  form a complete intersection of two cubics.
 These cubics formed by the unions of triples of lines in 
${{\bf P}^2}^*({\bf C})$ corresponding 
to triple of lines in ${{\bf P}^2}^*({\bf F}_3)$
passing through a point of $\ell$. 
Indeed, for a given $P_1,P_2 \in \ell$ and a point $Q$ on 
${{\bf P}^2}^*({\bf F}_3)$ outside of 
$\ell$, there are exactly 2 lines in ${{\bf P}^2}^* ({\bf F}_3)$ 
intersecting at this point and passing respectively through 
$P_1$ and $P_2$. The same incidence relation is valid on 
${{\bf P}^2}^*({\bf C})$.
\par In the case d) the homogeneous system has rank 9, i.e. 
the corresponding system does not define a torus.
\par Non essential tori correspond to subarrangements with number
of lines divisible by $m=3$ (cf. 2.6.3). 
There are 12 triples of lines
corresponding to each of triple points each defining a 2-torus.
A collection of 6 lines should have 4 triple points but the arrangement
of this example does not contain such subarrangements. 
\par Therefore we have 16 2-dimensional tori. In the abelian cover 
of ${\bf C}^2$ which sends each generator of $H_1({\bf C}^2-{\cal C})$ 
to a generator of ${\bf Z}/n{\bf Z}$ contributing tori are the 
essential torus of this arrangement and subtori corresponding to
subarrangements formed by triple of lines defined by the triple points.
Each torus contributes ${(n-1)(n-2)} \over 2$ to the Betti number i.e.
the total Betti number of this cover is  $16 \times 
{{(n-1)(n-2)}\over 2}$. 
These tori can be explicitly described as follows. Defining equations 
of non essential tori are products of 3 generators $t_i$'s corresponding
to a triple of points in ${\bf F}_3^2$ belonging to a line 
with the rest of $t_i$ is 1. Each of essential tori is given by 
9 equations $t_it_jt_k=1$ where $(i,j,k)$ are the triples of points 
${\bf F}_3^2$ (which interpreted as the lines of the arrangement)
which belong to lines not passing through a fixed point at infinity.
\par If the point at infinity is $(1,-1,0)$, then the lines 
not passing through it are: 
$x+z=0,x-z=0,x=0,y=0,x-y=0,x-y+z=0,x-y-z=0,y+z=0,y-z=0$ i.e. the 
corresponding torus satisfies:
$$t_{20}t_{21}t_{22}=1, t_{10}t_{11}t_{12}=1,
t_{00}t_{01}t_{02}=1,t_{00}t_{10}t_{20}=1,
t_{00}t_{11}t_{22}=1,$$ $$t_{01}t_{12}t_{20}=1,
t_{02}t_{10}t_{21}=1,t_{02}t_{12}t_{22}=1,
 t_{01}t_{11}t_{21}=1 \eqno (3.3.4)$$
 where the points of the complement to $z=0$ (i.e. the lines in 
${{\bf P}^2}^*({\bf C})$) are labeled as:
$$(0,0),(0,1),(0,2),(1,0),(1,1),(1,2),(2,0),(2,1),(2,2)$$ 
The corresponding torus  can be parameterized as
$$t_{00}=t,t_{01}=s,t_{02}=t^{-1}s^{-1},t_{10}=s,
t_{11}=t^{-1}s^{-1},t_{12}=t,$$
$$ t_{20}=t^{-1}s^{-1},t_{21}=t,t_{22}=s.$$
\par The equations for other essential tori, corresponding to
choices of the point at infinity  as respectively: $(1,1,0),
(1,0,0),(0,1,0)$ can be obtained from (3.3.4) by applying 
linear transformation to the indices which takes $(1,-1,0)$ 
to respective point. For $(x,y) \rightarrow (x,-y)$ 
which takes $(1,-1,0)$ to $(1,1,0)$ we obtain:
 $$t_{20}t_{22}t_{21}=t_{10}t_{12}t_{11}=t_{00}t_{02}t_{01}=
t_{00}t_{10}t_{20}=t_{00}t_{12}t_{21}=$$ 
 $$t_{02}t_{11}t_{20}=t_{01}t_{10}t_{22}=t_{01}t_{11}t_{21}=
t_{02}t_{12}t_{21}=1$$
For $(x,y) \rightarrow (x,x+y)$ which takes $(1,-1,0)$ to 
$(1,0,0)$ we obtain:
$$t_{22}t_{20}t_{21}=t_{11}t_{12}t_{10}=t_{00}t_{01}t_{02}=
t_{00}t_{11}t_{22}=t_{00}t_{12}t_{21}=$$
 $$t_{01}t_{10}t_{20}=t_{02}t_{11}t_{20}=t_{02}t_{10}t_{21}=
 t_{01}t_{12}t_{20}=1$$
For $(x,y) \rightarrow (x+y,y)$ which takes $(1,-1,0)$ to 
$(0,1,0)$ we have:
$$t_{20}t_{01}t_{12}=t_{10}t_{21}t_{02}=t_{00}t_{11}t_{22}=
t_{00}t_{10}t_{20}=t_{00}t_{21}t_{12}$$
$$t_{11}t_{01}t_{20}=t_{22}t_{10}t_{01}=
t_{22}t_{02}t_{12}=t_{11}t_{21}t_{01}.$$
\smallskip \par \noindent
{\it Example 4.} 
Let us consider the curve of degree 4 which has 
one ordinary point of multiplicity $4$. Faces of the polytopes of 
quasiadjunction are 
 $H_1: x_1+x_2+x_3+x_4=1$ (resp. $H_2: x_1+x_2+x_3+x_4=2$).
 The number of points
$({{i_1} \over n}, {{i_2} \over n},{{i_3} \over n}, {{i_4} \over n})$ 
on $H_1$ (resp. $H_2$) is ${(n-1)(n-2)(n-3)} \over 6$ (resp. 
 ${1 \over 3}(n-1)(2n^2-4n+3)$). The ideal corresponding to the 
polytope of quasiadjunction with the face $H_1$ (resp. $H_2$) is
${\cal M}^2$ (resp. $\cal M$ the maximal ideal of the local ring)
 and the level of  the supporting face
$H_1$ (resp. $H_2$) is 1 (resp. 2). Moreover $dim H^1 ({\bf P}^2, {\cal
J}_{{\cal M}^{3-l}} (4-3-l))$ is $2$ (resp. $1$) for $l=1$ 
(resp. for  $l=2$). Hence the irregularity of the cover corresponding to
homomorphism $H_1({\bf P}^2-\cup_{i=1,2,3,4} L_i) \rightarrow {\bf
Z}/n{\bf Z}$ is equal to 
$$2 \times {1 \over 6}(n-1)(n-2)(n-3)+{1 \over 3}(n-1)(2n^2-4n+3)+4
{1 \over 2}(n-1)(n-2)=$$ $$(n-1)(n^2-n-1). \eqno (3.3.5)$$  
This implies that the characteristic variety in this case is just 
 $$t_1t_2t_3t_4=1$$ 
The latter contains $(n-1)^3-(n-1)(n-2)=(n-1)(n^2-3n+3)$ points with 
coordinates in $\mu_n$ and the Betti number of the branched cover
from Sakuma's formula is equal to $2(n-1)(n^2-3n+3)+4(n-1)(n-2)=
2(n-1)(n^2-n-1)$.
\smallskip \par \noindent
{\it Example 5.} Let us consider the arrangement formed by 12 lines
which compose 4 degenerate fibers 
in a Hesse pencil of cubics formed by a non singular cubic curve 
and its Hessian. For example one can take the following
pencil: $$x^3+y^3+z^3-3\lambda xyz=0$$
This arrangement has 9 points of multiplicity 4 (inflection points
of non singular cubic). In ${{\bf C}^*}^{12}$ there are 10 tori of dimension
3 which are defined by 9 quadruples of lines corresponding to 
9 quadruple points and one 3-torus corresponding to 
the whole configuration. Contribution into the first Betti number 
an abelian cover also comes from 94 tori of dimension 2:
2-tori corresponding to triples of lines forming each of 9 quadruple
points (total $36$ 2-tori), 2-tori corresponding to  
configurations of 9 lines formed by triples of 4 special fibers of the 
pencil (total 4 2-tori) and 54 2-tori corresponding to configurations of
6 lines passing through 4 inflection points no three of which belong 
to a line (since 
the choice of 4 points must be made among points of affine space over
${\bf F}^3$ the ordered collection can be made in $9 \times 8 \times 6 
\times 3$ way and $54=9 \cdot 8 \cdot 6 \cdot 3 /24$). In particular, the 
irregularity of the cover with the Galois group $({\bf Z}/3 {\bf Z})^2$
is equal to 154 (cf. [I]).
Indeed the contribution of each 2-torus into the first Betti number is
2 and in the case of 3-tori the contribution is 6,
since the 3-torus contains 6 points with coordinates $i \over 3$.
Since the depth of 3-tori is 2 the first Betti number is equal to 
$6 \times 10 \times 2+94 \times 2$.
\section{The structure of characteristic varieties of algebraic
curves.}

In this section we describe sufficient conditions for the vanishing of 
cohomology of linear systems which appear in description of 
characteristic varieties given in section 3. This, therefore, 
yields conditions for absence of essential components.
In the cyclic case one obtains triviality of Alexander polynomial.
\subsection{Absence of characteristic varieties for curves with 
small number of singularities.}
\subsubsection{}
\bigskip 
\begin{theo}
Let ${\cal C}$ be a plane curve as above. 
Suppose that $\rho: Y \rightarrow  {\bf P}^2$ is obtained by a 
sequence of blow ups such that 
the proper preimage $\tilde {\cal C}$ of ${\cal C}$ in $Y$ 
has only normal intersection with the exceptional set 
and satisfies $\tilde {\cal C}^2 >0$. Then ${\cal C}$ has no essential 
characteristic subvarieties.   
\end{theo}
\smallskip  \par \noindent
\begin{coro} 1. Let $\cal C$ be an irreducible 
curve which has ordinary cusps and nodes as the only singularities.
If the number of cusps is less than ${{d^2} \over 6}$ then the 
Alexander polynomial of $\cal C$ is equal to 1. 
\par 2. Let $\cal H$ be an arrangement consisting of $d$ lines and
which has $N$ points of multiplicity $m$. Let $l(\delta)$ be the 
level of a face of quasiadjunction for the complement to $\cal H$. 
If $d^2 >m^2 N$ then the superabundance is zero for 
the system of curves 
of degree $d-3-l(\delta)$ which local equations belong to 
the ideal of quasiadjunction corresponding to $\delta$ at the points 
of multiplicity $m$.
\end{coro}
\bigskip \par \noindent
{\it Remark.} One can compare corollary 1 
with Nori's results (cf. [N]). The latter yields that the fundamental 
group of the complement to a curve of degree $d$ with $\delta$ nodes
and $\kappa$ cusps is abelian if $d^2 > 6 \kappa +2 \delta$ while 
a weaker inequality $d^2 > 6 \kappa$ yields the triviality of 
the Alexander polynomial. For example, for the branching curve of 
a generic projection  of a smooth surface of degree $N$ in 
${\bf P}^3$ one has $d^2 > 6 \kappa$ for $N >4$ but $d^2 < 6 \kappa
+2 \delta$ for $N >2$.  The fundamental groups of these curves 
are non abelian for $N >2$ and the Alexander polynomial 
for $N=3,4$ is equal to $t^2-t+1$ (cf. [L2]). 
\bigskip \par \noindent
{\bf Proof of the theorem.} 
We should show that for any contributing
face of quasiadjunction $\delta$ we have 
 $dimH^1({\bf P}^2, {\cal A}_{\delta}(\Sigma d_i-3-l({\delta})))=0$.
If $\rho: Y \rightarrow {\bf P}^2$ is a  blow up of  
${\bf P}^2$, satisfying conditions of the theorem, 
then we have: 
$${\cal A}_{\delta}(\Sigma_i d_i -3-l({\delta}))=
\rho_*(\omega_Y \otimes {\cal
O}_Y(\gamma \tilde {\cal C})\otimes {\cal O}_Y(\Sigma \epsilon_k E_k)) \eqno (4.1.1)$$
for some rational $\gamma>0$ and 
$0 \le \epsilon <1$. More precisely, 
$\gamma={1 \over {\Sigma_i d_i}} \cdot \Sigma_i d_i(1-{{j_i+1} \over {m_i}})$ 
for some choice of $(...,{{j_i+1} \over m_i},...)$ belonging to the 
face  $\delta$
(with $\epsilon_k$ a priori depending on this choice). Indeed, 
from (2.3.6) and the discussion after, 
the multiplicity $f_k(\phi)$ along an exceptional curve 
$E_k$ of the pull back on 
$Y$ of a germ 
in the ideal of quasiadjunction with parameters 
$(j_1,..,j_r \vert m_1,...,m_r)$ such that $({{j_1+1}\over m_1},...,
{{j_r+1} \over m_r}) \in \delta$ satisfies:
   $$f_k(\phi) \ge [\Sigma_i a_{k,i} (1-{{j_i+1} \over m_i})-c_k] 
\eqno (4.1.1.1)$$
Hence $$A_{\delta}=\rho_* (\otimes_k  
{\cal O}_Y( (c_k-[\Sigma_i a_{k,i}-a_{k,i}({{j_i+1} \over m_i})]) E_k)) \eqno
(4.1.1.2)$$
We have: $l(\delta)=\Sigma_i d_i {{j_i+1} \over {m_i}}$ and 
${\cal O}_{{\bf P}^2}(C_i)={\cal O}_{{\bf P}^2}(d_i)$.
 Therefore 
$$ A_{\delta}(\Sigma_i d_i-3-l(\delta))=$$ $$\rho_* (\otimes_k  
{\cal O}_Y( (c_k-[\Sigma_i a_{k,i}-a_{k,i}({{j_i+1} \over m_i})])
E_k)) \otimes {\cal O}_{{\bf P}^2}(-3) 
\otimes {\cal O}_{{\bf P}^2}(\Sigma_iC_i(1-{{j_i+1} \over m_i})) 
\eqno (4.1.1.3)$$
Since $$ \omega_Y= \otimes_k {\cal O}_Y(c_kE_k) \otimes
\rho^*({\cal O}_{{\bf P}^2}(-3))$$  and
$$\otimes_k {\cal O}_Y( a_{k,i}(1-{{j_i+1} \over m_i})E_k)\otimes 
{\cal O}({\tilde {\cal C}})^{{{d_i}  \over {\Sigma d_i}}}=\rho^*(
 {\cal O}_{{\bf P}^2}(C_i(1-{{j_i+1} \over m_i}))) \eqno (4.1.1.4)$$
(because  ${\cal O}_{{\bf P}^2}(C_i)={\cal O}_{{\bf P}^2}(1)^{d_i}=
{\cal O}_{{\bf P}^2}({\cal C})^{{d_i} \over {\Sigma_i d_i}}$),
 we see that (4.1.1.3) yields (4.1.1) 
with $\epsilon_k=\{\Sigma_i a_{k,i} (1-{{j_i+1} \over m_i}) \}$ 
where $\{ x \}=x-[x]$ is the fractional part.
\par The Kawamata-Viehweg vanishing theorem (cf. for example [Ko]) 
implies that the 
cohomology of the sheaf $\omega_Y \otimes 
{\cal O}_Y(\gamma \tilde {\cal C}) \otimes 
{\cal O}_Y(\Sigma \epsilon_k E_k)$ is trivial in positive dimensions
 if $\tilde {\cal C}$ is big 
and nef. But this  follows from the assumptions of the theorem.
Finally, the exact sequence $0 \rightarrow E_2^{1,0} \rightarrow 
H^1(Y,{\cal F})$ of lower degree terms in the 
Leray spectral sequence $H^p({\bf P}^2,R^q \rho_*{\cal F})
\Rightarrow H^{p+q}(Y,{\cal F})$ for the sheaf ${\cal F}=
\omega_Y \otimes 
{\cal O}_Y(\gamma \tilde {\cal C}) \otimes 
{\cal O}_Y(\Sigma \epsilon_k E_k)$
yields $H^1({\bf P}^2, {\cal A}_{\delta}(\Sigma_i d_i-3-l(\delta)))=0$.

\bigskip \noindent
{\bf Proof of the Corollary.} 
For each blow up at an ordinary  point of multiplicity $m$ of a curve
$C$ we have 
${\tilde C^2}=C^2-m^2$ where $\tilde C$ is the proper preimage of $C$. 
Hence ${\cal O}(\tilde {\cal H})$ is big if $d^2>m^2N$. 
The case of ordinary cusps is similar. 
\subsubsection{Dimensions of components of characteristic 
varieties.}
\bigskip
This theorem imposes restrictions on the  
dimensions of the contributing faces and hence on the dimensions
of characteristic varieties. For example, let us consider an arrangement of
lines with at most triple points as singularities. Then each contributing
face is the intersection of hyperplanes defining
the only local polytope of an ordinary triple point.
These hyperplanes are given by the equations of the form 
$x_i+x_j+x_k=1$ where $(i,j,k)$ 
are the indices corresponding to the lines through the triple point.
The matrix of this system therefore has the property that  
in each row only 3 non zero entries are equal to 1, any two rows have
at most one non zero entry in the same column and the number of rows is
at least ${d^2} \over 9$ (since by the corollary only in this case  one can
get a contributing face with $H^1 \ne 0$). 
In particular, the number of non zero entries in the matrix is at 
least ${{d^2} \over 3}$. 
The rank of this system is at
least $d \over 3$. Indeed, the matrix contains a column 
with at least $d \over 3$ 1's. On the other hand, if a column $s$ 
contains $k$ non zero entries in rows $v_1,..,v_k$, then these rows 
are linearly independent since the left hand side of 
a relation $\Sigma \lambda_i v_i=0$ 
has $\lambda_i$ as the entry of the column different from the $s$-th 
and in which $v_i$ has a non zero entry.
In particular in the arrangement of $d$ lines, the dimension of
characteristic variety is at most ${2d} \over 3$. For the components
containing a trivial character, the dimension is at most 2, since 
by [Ar] such component induces the map onto ${\bf P}^1$ minus three 
point such that pull backs of rank one local systems from the latter
form the component  (cf. footnote in section 1 and [LY]).
\smallskip \par \noindent
{\bf Problem.} Let $m$ be the maximal number of components which
meet at a singular point of $\cal C$. Is it true that the dimension of each
component of the characteristic variety of ${\cal C}$ is at most $m-1$?
\bigskip \par \noindent
\subsection{Translations of the tori forming characteristic
varieties and the degrees of irreducible components of $\cal C$} 
\smallskip \par \noindent
\begin{theo} Let $a_{\delta}$ 
be g.c.d. of non zero minors of maximal order in 
a system of equations with integer coefficients which set of solutions 
contains a face of quasiadjunction $\delta$ with the dimension 
of the set of solutions equal to $\dim \delta$. Let $d_i=deg C_i$. 
\par \noindent
a) Each irreducible essential component of $Char
{\cal C}$, belongs to a coset of a subtorus of 
$\Hom (H_1({\bf C}^2-{\cal C}), {\bf C}^*)$ 
of order dividing $a_\delta$.
In particular each irreducible essential component belongs to 
a coset of order equal to the order of the face (cf. 2.4.2 and 
2.6.2).
\par \noindent 
b) Each essential component belongs to a codimension 
one subtorus of $\Hom (H_1({\bf C}^2-{\cal C}), {\bf C}^*)$
translated by a point of the order dividing 
$\varrho=g.c.d. (d_1,...,d_r)$.
\end{theo}
\par \noindent 
{\bf Proof.} Any component of a characteristic variety 
is a Zariski closure of the image under the exponential map 
of a contributing face $\delta$ with the equation
 $d_1x_1+...+d_rx_r=l(\delta)$. This yields b) (it also 
follows from Cor.3.3 from [L3]).
\par It follows from theorem 3.1 c) that irreducible essential 
component is a connected component of the subgroup of 
 $\Hom (H_1({\bf C}^2-{\cal C}), {\bf C}^*)$ belongs to a subgroup  
$$\chi_1=...=\chi_r=1 \eqno (4.2.1.1)$$ 
where $\chi_i$ is the character of 
$\Hom (H_1({\bf C}^2-{\cal C}), {\bf C}^*)$ having form 
$exp (L_i)$ where $L_i$ is the form with integer 
coefficients such that $L_i=l_i, l_i \in {\bf Z}$ are the equations defining 
the face of quasiadjunction. The order of the group of cosets 
of the group (4.2.1.1) by its connected component of identity
is the order of the torsion of its group of characters.
The latter is $\Char(\Hom (H_1({\bf C}^2-{\cal C}), {\bf
C}^*))/(\chi_1,... \chi_r)$. This yields a).
\addtocounter{subsubsection}{1}
\subsubsection{Linear systems corresponding to different faces of 
quasiadjunction.}
Another byproduct of results in section 3 is equality of 
superabundances of linear systems of curves defined by 
rather different local conditions.
\bigskip
\begin{prop} a) Let $\delta$ and $\delta'$ be two 
faces of global polytopes quasiadjunction such that the Zariski closures of 
$exp(\delta)$ and $exp(\delta')$ coincide. 
Then if $\delta$ is a contributing face then $\delta'$ is 
also contributing and 
$H^1({\cal A}_{\delta}(d-3-l(\delta))=H^1({\cal A}_{{\delta}'}
(d-3-l({\delta}'))$. 
\par \noindent b) Let $\alpha \in {\bf Q}$ is such that 
$\alpha \cdot g.c.d (d_1,...,d_r)$ is the level 
of a face of  quasiadjunction $\delta$ and 
$\sigma \in Gal({\bf Q}(exp(2 \pi i \alpha)) /{\bf Q})$ 
such that $\sigma(exp (2 \pi i \alpha))=exp(2 \pi i \beta)$
 with $0 <\beta <1$. Then 
$\beta$ is equal to 
 ${{l({\delta}')} \over {g.c.d (d_1,..,d_r)}}$
for some face of quasiadjunction $\delta '$
and $H^1({\cal A}_{\delta}(\Sigma_i d_i-3-
l(\delta)))=H^1({\cal A}_{{\delta}'}(d-3-l({\delta}')))$. 
\end{prop}
\smallskip \par \noindent
{\bf Proof.} a) Since Zariski closures of $exp(\delta)$ and
$exp(\delta')$ are the same the corresponding to $\delta$ and 
$\delta'$ components of the 
characteristic variety are the same. If the depth of this
component of  characteristic variety is
 $i$, then the dimension of each of the cohomology group in the 
statement equals to  $i$ and the result follows. 
\par \noindent 
b) Since i-th characteristic variety is defined over $\bf Z$ the 
Galois group $Gal(\bar {\bf Q}/{\bf Q})$ acts on the set
of its irreducible components. Since  
irreducible component corresponding to $\delta$ is a translation 
by $exp (2 \pi \sqrt {-1} \alpha)$
of a subgroup of $H^1({\bf C}^2-{\cal C},{\bf C}^*)$
defined over ${\bf Z}$, it follows that  
$\sigma$ takes the component corresponding to $\delta$ 
into translation of the same subgroup 
by $exp (2 \pi \sqrt {-1} \beta)$.
This translation is a Zariski closure of 
$exp(2\pi \sqrt {-1} \delta ')$ for some face $\delta '$. 
It does satisfy the conclusions of b).
\smallskip \par \noindent
\subsubsection{Remarks.} 
\bigskip \par \noindent 1. For irreducible curves 
the order of each face of quasiadjunction is a root of 
a local Alexander polynomial. So the divisibility theorem 
from [L1] is a special case of 4.2.1.
\par \noindent  2. One of the consequences of 4.2.1 is a 
non trivial restriction on an abstract group which is necessary to satisfy in 
order that the group can be realized as the fundamental group of 
an arrangement. For example $t_1 \cdot \cdot \cdot t_r=-1$ 
cannot be a component of 
characteristic variety of arrangement of $r$ lines since it is 
cannot belong to an intersection of subgroups of ${{\bf C}^*}^r$.
\par \noindent 3. As an illustration to 4.2.2b), 
let us consider an irreducible 
curve of degree $d$ with singularities locally isomorphic to singularity 
$x^2=y^5$.
 If the linear system consisting of curves of 
degree $d-3-{d \over {10}}$ with local equations belonging to ideals of
quasiadjunction of all singular points 
corresponding to the constant of quasiadjunction $1
\over {10}$ is superabundant, then the linear system of curves of degree
$d-3-{{3d} \over {10}}$ with local equations in the ideals of
quasiadjunction corresponding to $3 \over {10}$ is also superabundant
and the superabundances are equal.
\par 4. An example of faces of quasiadjunction with the same 
Zariski closure of the images of the exponential map 
as in a) of the proposition is given 
by  $x_1+...x_m=i, x_1+..+x_m=j, 0< i,j \le m-2$ which are the faces of
quasiadjunction for the complement to $m$ lines through a point
(cf. 2.5 example 3).    
\section{Resonance conditions for rank one local systems  on 
complements to line arrangements.}
\subsection{Complexes associated with arrangements.}
Let ${\cal L}=\cup_{i=1,...,r} L_i$ be an 
arrangement of lines in ${\bf C}^2$. We shall assume for convenience 
(cf. (1.2.3)) that
the line at infinity is transversal to all lines in ${\cal L}$. 
Let $l_i(x,y)=0$ be the equation of $L_i$ and 
$\eta_i={1 \over {2 \pi i}} {{dl_i} \over l_i}$. Let $A^i (i=0,1,2)$
be the subspace generated by the forms 
$\eta_{j_1} \wedge ... \wedge \eta_{j_i}$ in the space 
$H^0({\Omega}^i(*{\cal L}))$ of meromorphic forms with poles along ${\cal
L}$. Let $$\omega=\Sigma \eta_i \cdot s_i, s_i \in {\bf C}
\eqno(5.1.1)$$ The exterior product with $\omega$ defines the complex:
 $$A^{\bullet}_{\omega}: 0 \rightarrow A^0 \rightarrow A^1 \rightarrow A^2
\rightarrow 0\eqno (5.1.2)$$ 
If $s_i=0, (i=1,...,r)$, then the cohomology groups of $A^{\bullet}$ 
are isomorphic to the cohomology groups of
${\bf C}^2-{\cal L}$ (cf. [Br]). On the other hand,
the collection $s_{*}=(s_1,...,s_r)$ defines the map $\pi_1({\bf C}^2-{\cal L})
\rightarrow {\bf C}^*$ which sends $\gamma_i$ (cf. (1.1)) to 
$exp (2 \pi \sqrt {-1} s_i)$ and hence the local system which we shall 
denote ${\cal A}_{s_*}$. A theorem from  [ESV] (p.558), in the case of
line arrangements, asserts that  
 $$H^i({\cal A}_{s_*})=H^i(A^{\bullet}_{\omega}) \eqno (5.1.3)$$
provided the following non resonance condition is satisfied. For any
point singular point $P$ of $\cal L$ of multiplicity $m > 2$, 
if the lines through $P$ are $l_{i_1},...,l_{i_m}$, then  
$$s_{i_1}+...+s_{i_m} \ne n \in {\bf N}-{0} \eqno (5.1.4)$$
\begin{theo} The isomorphism (5.1.3) takes 
place, provided $$(exp(2 \pi i s_1),..., exp (2 \pi i s_r)) \eqno (5.1.5)$$ 
does not belong to the characteristic variety $Char_1$ of ${\bf C}^2-{\cal L}$.
\end{theo}
\subsection{Remarks.} 1. It is easy to construct examples of local
systems for which (5.1.4) is violated but for which (1.5.3) takes place.
Indeed, the image under the exponential map onto the
torus ${{\bf C}^*}^r$ of those $(s_1,...,s_r)$ which violate (5.1.4) 
is a union of codimension 1 tori in ${{\bf C}^*}^r$. On the other hand, 
the characteristic varieties typically have rather small dimension 
relative to $r$ (cf. (3.3) and (4.1)).
\par \noindent 2. For arrangements of arbitrary dimension in ${\bf C}^n$ 
(with $l_i$ denoting 
the equations of hyperplanes of the arrangement, rather than lines)
we have  $$H^1(A^{\bullet}({\cal L}),\omega)=0, \eqno (5.2.1)$$ 
provided condition (5.1.5) of the theorem is met. Indeed, 
for a generic plane $H \subset {\bf C}^n$ the map 
$A^i({\cal L}) \rightarrow A^i({\cal L} \cap H)$ 
induced by inclusion is isomorphism for $i=0,1$ and injective for $i=2$.
The latter follows from the Lefschetz theorem since $A^i({\cal L})=H^i({\bf
C}^n-{\cal L})$ by [Br]. This yields the isomorphism
of cohomology of the complexes $A^{\bullet}$ for ${\cal L}$ and 
${\cal L} \cap H$ and hence (5.2.1).
Therefore the theorem from (5.1) (and also from (5.4))
holds for arbitrary arrangements. 
\subsection{Proof} 
We will derive this theorem from 
the following:
\par \noindent
a)$H^0(A^{\bullet}_{\omega})=H^1(A^{\bullet}_{\omega})=0$
\par \noindent 
b)The euler characteristics of both $A^{\bullet}_{\omega}$ 
and ${\cal A}_{s_*}$   
are equal to $e({\bf C}^2-{\cal L})$.
\par \noindent
To show b) note that the euler characteristic of $A^{\bullet}$ is $e({\bf C}^2-{\cal L})$
by [Br] and for $H^i({\cal A}_{s_*})$ this can be seen by looking at 
the cochain complex of ${\cal A}_{s_*}$, i.e. (here $\chi_{s_*}$ is the
character of the fundamental group defining the local system ${\cal A}$)
 $$C^i((\widetilde
{{\bf C}^2-{\cal L})} \otimes_{\chi_{s_*}} {\bf C})={\bf C}^{b_i} 
\otimes {{\bf C}[H_1({\bf C}^2-{\cal L})]} \otimes_{\chi_{s_*}} {\bf C}
 \eqno (5.3.1)$$
since the multiplicity of a representation in the regular representation
of an abelian group is 1.
\par On the other hand a) follows since
 $dim H^1(A^{\bullet},\omega) \le dim H^1({\cal A}_{s_*})$ (cf. [LY]) 
 and the latter  group is trivial if (5.1.5) is satisfied.
\par Now the theorem 5.1 is a consequence of  a),b) and (1.4.1.1).
\smallskip \par \noindent
\subsection{Complexes (5.1.2) with non vanishing cohomology.}
Complexes $A^{\bullet}$ for a fixed arrangement ${\cal L}$ are parameterized
by the space $H^0({\bf P}^2,\Omega^1(log({\cal L} \cup L_{\infty}))
=H^1({\bf C}^2-{\cal L},{\bf C})$. Let 
$${\cal V}_i =\{ \omega \in H^0({\bf P}^2,\Omega^1(log({\cal L} 
\cup L_{\infty})))
\vert H^1(A^{\bullet},\omega) \ge i \} \eqno (5.4.1)$$
\par \noindent
\begin{theo} 
${\cal V}_i$ is a union of linear space of dimension 
$i+1$. There is the one to one correspondence between these linear spaces and 
irreducible components of the characteristic variety $V_i$ 
containing the trivial character which have 
a positive dimension. For such an irreducible component  
of ${\cal V}_i$ the map which assigns to a form $\omega$ the 
point $(...,exp(2 \pi i Res_{L_i} \omega),...)$ in 
$H^1(\pi_1({\bf C}^2-{\cal L}),{\bf C}^*)$ (using identification (1.4.1))
is the universal cover of  the corresponding component of $V_i$.
\end{theo}
\smallskip \par \noindent 
{\bf Proof.} If $dim H^1(A^{\bullet},\omega)>0$, i.e.
there exist linearly independent with $\omega$ form $\eta$ such that 
$\eta \wedge \omega=0$, then for any form $\omega'$ in the space 
spanned by $\omega$ and $\eta$ one has $dim H^1(A^{\bullet},\omega')>0$.
Let ${\cal V}$ be an irreducible component of ${\cal V}_i$ in (5.4.1) 
containing $\omega$ and having dimension $k \ge 2$. 
For the local system $L_{\omega'}$ corresponding to each  $\omega'
\in {\cal V}$ we 
have $dim H^1(L_{\omega'})>0$. Indeed we can assume   
that $\omega'$ is  generic since this 
only decrease $H^1(L_{\omega'})$ (cf. [LY]). On the other hand for 
generic $\omega'$, according to [ESV], we have 
$dimH^1(A^{\bullet}, \omega')=dim H^1(L_{\omega'})$. 
Therefore $L_{\omega}$ belongs to an irreducible component, say $V$, of the 
characteristic variety of ${\bf C}^2-{\cal L}$. 
Since the exponential map is a local homemorphism this component has
the dimension equal to at least $k$. In fact the dimension of this 
component is exactly $k$. Assume to the contrary that this dimension is  
$l >k$  and let $f: {\bf C}^2-{\cal L}
 \rightarrow {\bf P}^1-\cup_{i=1}^{i=k+1} p_i$
be the map on a curve of general type (cf. (1.4.2) and [Ar], Prop. 1.7)
corresponding to the component $V$. Then the pull back
of form $H^0({\bf P}^1, \Omega^1(log(\cup p_i)))$ gives $l$-dimensional
space of forms on ${\bf C}^2-{\cal L}$ for which the wedge with $\omega$
is zero (note that the map $f^*$ is injective on $H^1$) and we have a 
contradiction.
Let $t_1^{a_{1,j}}...t_r^{a_{r,j}}=1 (j=1,...s)$ be the equations defining 
$V$ (cf. 4.2). Then $\omega$ belongs to the union of affine subspaces 
of $H^0({\bf P}^2,\Omega^1(log ({\cal L} \cup L_{\infty}))$
 given by $Q_j=\Sigma a_{i,j}x_i=n_j, n_j \in {\bf Z}, j=1,...,s$.
Since $dim H^1(A^{\bullet},\lambda \omega)=
dim H^1(A^{\bullet}, \omega), \lambda \in {\bf C}^*$ we see that $n_j=0$
for any $j$. Hence ${\cal V}$ is a linear space of 
dimension $k$ (and $i=dim H^1(A^{\bullet},\omega)=dimH^1(L_{\omega})=k-1$).    
\smallskip \par \noindent
\subsection{ Combinatorial calculation of characteristic varieties.}
\par A consequence of the theorem 5.4 is that the irreducible components
of the characteristic varieties containing the identity element 
of $H^1({\bf C}^2-{\cal L},{\bf C}^*)$
 are determined by the cohomology of the complex (5.1.2). 
$H^1(A^{\bullet},\omega)$ is the quotient of $\{\eta \in A^1 \vert 
\eta \wedge \omega =0 \}$ by the subspace spanned by $\omega$ and 
can be calculated as follows. 
\par It is easy to see that a 2-form is cohomologous to zero 
iff its integrals over all 2-cycles belonging to small
balls about the multiple points of the arrangement are zeros. 
The group of such 2-cycles near a point which is the intersection of 
the lines $l_{i_1},...,l_{i_m}$ are generated by 
$\gamma_{i_j} \times ({\gamma_{i_1}+....+\gamma_{i_m}})$. 
If $(\Sigma A_i \eta_i) \wedge \omega$ is cohomologous to
zero in $\Omega^2({\bf C}^2-{\cal L})$ 
then vanishing of  $\int A_j \eta_j \wedge 
s_i \eta_i$ over those 2-cycles yields: $$A_j(\Sigma
s_i)-(\Sigma A_j)s_j=0 \eqno (5.5.1)$$ 
Therefore we obtain $$A_j=C_{\upsilon}s_j (if \Sigma_{\upsilon \in l_j}
s_j \ne 0) \eqno (5.5.2)$$
 $$\Sigma A_j=0  (if \Sigma_{\upsilon \in l_j}
s_j=0) \eqno (5.5.3)$$ for vertices  $\upsilon$ of the arrangements.
If we are looking for essential components of the characteristic 
variety (which we always can assume) then $s_i\ne 0$ and condition 
 (5.5.2) can be replaced by 
    $${{A_j} \over {s_j}}= {{A_{j'}} \over {s_{j'}}}
\ \ \ if \ \ \Sigma_{\upsilon \in l_j} s_j \ne 0 \eqno (5.5.4)$$
Now for each subset of the set of vertices such that the 
system of equation (5.5.3) and supplementing it by equations (5.5.4)
 for vertices outside of selected subset has a solution non proportional
to $(s_1,..,s_r)$ we obtain a component $\cal V$ and hence corresponding 
component of the characteristic variety.
We leave as an exercise to the reader to work out calculations of the 
characteristic varieties for the examples from section 3 using 
this method.



\begin{thebibliography}
\bigskip
\bigskip 
\bibitem*[Ar] Arapura, D., (1997) Geometry of cohomology support loci for local
systems, I, {\it J. Alg. geom}, {\bf 6}, 563--597.
\bibitem*[BE] Buchsbaum D. and Eisenbud D., (1977) 
What annihilates a module?,
{\it J.Algebra}, {\bf 47}, 231-243.
\bibitem*
[BHH] Barthel, G., Hirzebruch, F, Hofer T., (1987) 
{\it Geradenkonfiguarationen und  Algebraishe fl\"achen}, 
 Vieweg Publishing, Wiesbaden.

\bibitem*
[BL] Blass P. and Lipman J., (1979) Remarks on adjoints and arithmetic genus
of algebraic varieties, {\it Amer. J. of Math}, {\bf 101},  331--336.
\bibitem*
[Br] Brieskorn E., (1973) Sur les groups des tresses, In Sem. Bourbaki,
{\it  Lect. Notes in Math.} {\bf 317}  Berlin, Heidelberg,
New York, Springer, 21-44. 
\bibitem* 
[CE] Cartan H. and Eilenberg S., (1956) {\it Homological Algebra}, 
Princeton University Press, Princeton, N.J. 
\bibitem*
[CS1] Cohen D. and Suciu A., (1995) On Milnor fibrations of arrangements, 
{\it J.London Math. Soc}, {\bf 51}, 105-119. 
\bibitem*
[CS2] Cohen D. and Suciu A., (1999) Characteristic varieties of 
arrangements. {\it Math. Proc. Cambridge Phil. Soc.}
{\bf 127}, 33--53. 
\bibitem*
[De] Deligne P., (1971) 
Theorie de Hodge II, {\it Publ. Math. I.H.E.S. } {\bf 40}, 5-58.
\bibitem*
[ESV] Esnault H., Schechtman V., and Viehweg E., (1992) Cohomology of local
systems on the complement to hyperplanes. 
{\it Inv. math.} {\bf 109}, 557-561. 
\bibitem*
[F] Falk M., (1997) 
Arrangements and cohomology. {\it Ann. Comb.}, {\bf 1},  135--157.
\bibitem*
[GR] Grauert H. and Riemenschneider O., (1970) {\it Inv. Math.} 
{\bf 11}, 263-292.
\bibitem*
[Ha] Hartshorne R., (1977) {\it Algebraic Geometry}, Springer Verlag, 
Berlin Heidelberg New York.  
\bibitem*
[Hil] Hillman J., (1981) Alexander Ideals of Links {\it 
Lecture Notes in Math.} 
{\bf 895}. Springer Verlag, Berlin Heidelberg New York.
\bibitem*
[Hi] Hironaka E., (1997) Multi-polynomial invariants of plane algebraic
curves, {\it Singularities and Complex Geometry, Studies in Advanced
Mathematics}, {\bf 5}, AMS and International Press, 67-74.
\bibitem*
[I] Ishida M., (1983) The irregularity of Hirzebruch's examples of surfaces
of general type with $c_1^2=3c_2$. {\it Math. Ann.} {\bf  262}, 407-420.
\bibitem*
[vK] van Kampen E.R., (1933) On the fundamental group of an algebraic 
curve, {\it Amer. J. of Math}, {\bf 55}, 255-260.
\bibitem*
[Ko] Kollar J., (1995) Shafarevich maps and automorphic forms, Princeton
University Press.
\bibitem*
[L1] Libgober A., (1982) Alexander polynomials of plane algebraic curves
and cyclic multiple planes, {\it Duke Math. J,} {\bf 49}. 833--851.
\bibitem*
[L2] Libgober A., (1983) Alexander invariants of plane algebraic curves, 
{\it Proc. Symp. Pure Math.} {\bf 40}, AMS. Providence, RI, 135-143. 
\bibitem*
[L3] Libgober A., (1992) On homology of finite abelian coverings, 
{\it Topology and Applications}, {\bf 43}, 157-166.
\bibitem*
[L4] Libgober A., (1994) Groups which cannot be realized as fundamental
groups of the complements to hypersurfaces in ${\bf C}^N$, {\it Algebraic
Geometry and Application, C.Bajaj editor,}  Springer Verlag, Berlin
Heidelberg New York, 203-207.
\bibitem*
[L5] Libgober A., (1996) Position of singularities of hypersurfaces and the
topology of their complements, {\it J.Math. Sciences}, {\bf 82}, 3194-3210. 
\bibitem*
[L6] Libgober A., (1999) Abelian Covers of projective plane, 
{\it Singularity
theory (Liverpool, 1996), xxi, 
London Math. Soc. Lecture Note Ser.,} {\bf 263},
Cambridge Univ. Press, Cambridge, 281-289.
\bibitem*
[L7] Libgober A., (2001) Hodge decomposition of Alexander invariants.
Preprint, 2001.
\bibitem*
[LY] Libgober A., Yuzvinsky S., (2000) Cohomology of the
Brieskorn-Orlik-Solomon algebra and local systems, 
{\it Arrangements---Tokyo 1998, Adv. Stud. Pure Math.}, {\bf 27}, 
Kinokuniya, Tokyo, 169-184. 
\bibitem*
[LV] Loeser F.,and Vaquie M., (1990) Le polynome d'Alexander d'une courbe
plane projective, {\it Topology} {\bf 29}, 163-173. 
\bibitem*
[N] Nori M., (1983) Zariski's conjecture and related problems. 
{\it Ann. Sci. Ecole Norm. Sup.} {\bf 16}, 30 5-344.
\bibitem*
[Sab] Sabbah C., (1990) Module d'Alexander et ${\cal D}$-modules, 
{\it Duke Math. Journal}, {\bf 60}, 729-814.
\bibitem*
[Sa] Sakuma M., (1995) Homology of abelian coverings of links and spatial
graphs, {\it Canadian Journal of Mathematics}, {\bf 17}, 201-224.
\bibitem*
[Se] Serre J.P., (1975) Algebre locale. Multiplicites, {\it Lecture Notes in
Mathemematics} {\bf 11}, (3rd edition), Springer Verlag, Berlin Heidelberg New York. 
\bibitem*
 [St] Steenrod N., (1945) Homology with local coefficients, {\it Ann. Math.},
 {\bf 44}, 610-627.
\bibitem*
[SW] Sumners D. and  Woods J., (1977) The monodromy of reducible plane
curves, {\it Inv. Math}. {\bf 40}, 107-141.
\bibitem*
[Z] Zariski O., (1971) {\it Algebraic surfaces}, Chapter.8, 
Springer Verlag, Berlin Heidelberg New York.
\bibitem*
[Zu] Zuo K., (1989) Kummer Oberlagerungen algebraischer Fl\"achen. {\it Bonner
Mathematische Schriften}, Bonn.

\end{thebibliography}
\end{document}